\newif\ifSIAM
\SIAMfalse

\ifSIAM
\documentclass[review]{siamart250211}
\else
\documentclass[letterpaper]{article}
\usepackage[margin=1.3in]{geometry}
\usepackage{amsthm}
\fi
\let\citet\cite

\usepackage{import}
\usepackage{amsmath}
\usepackage{mathtools}
\usepackage{amssymb}
\usepackage{mathrsfs}
\usepackage{xparse}
\usepackage{import}
\usepackage[normalem]{ulem}

\makeatletter
\@ifpackageloaded{unicode-math}{%
	\newcommand{\mymathbold}{\symbf}%
}{%
	\usepackage{bm}%
	\newcommand{\mymathbold}{\bm}%
}
\makeatother

\newcommand{\scrbar}[1]{\overline{\mathcal{#1}}}
\newcommand{\scrhat}[1]{\widehat{\mathcal{#1}}}
\newcommand{\scrtl}[1]{\widetilde{\mathcal{#1}}}

\ExplSyntaxOn

\cs_new_protected:Nn \bamboo_define:nnnnnN
{
	\cs_new_protected:cpx { #3 #1 #4 } { \exp_not:N #5{#6{#2}} }
}

\int_step_inline:nnn { `A } { `Z }
{
	\bamboo_define:nnnnnN
	{ \char_generate:nn { #1 } { 11 } }
	{ \char_generate:nn { #1 } { 11 } }
	{ bl }
	{ }
	{ \mymathbold }
	\use:n
	
	\bamboo_define:nnnnnN
	{ \char_generate:nn { #1 } { 11 } }
	{ \char_generate:nn { #1 } { 11 } }
	{ scr }
	{ }
	{ \mathcal }
	\use:n
	
	\bamboo_define:nnnnnN
	{ \char_generate:nn { #1 } { 11 } }
	{ \char_generate:nn { #1 } { 11 } }
	{ }
	{ hat }
	{ \widehat }
	\use:n
	
	\bamboo_define:nnnnnN
	{ \char_generate:nn { #1 } { 11 } }
	{ \char_generate:nn { #1 } { 11 } }
	{ scr }
	{ hat }
	{ \scrhat }
	\use:n
	
	\bamboo_define:nnnnnN
	{ \char_generate:nn { #1 } { 11 } }
	{ \char_generate:nn { #1 } { 11 } }
	{ }
	{ br }
	{ \overline }
	\use:n
	
	\bamboo_define:nnnnnN
	{ \char_generate:nn { #1 } { 11 } }
	{ \char_generate:nn { #1 } { 11 } }
	{ }
	{ tl }
	{ \widetilde }
	\use:n
	
	\bamboo_define:nnnnnN
	{ \char_generate:nn { #1 } { 11 } }
	{ \char_generate:nn { #1 } { 11 } }
	{ scr }
	{ br }
	{ \scrbar }
	\use:n
	
	\bamboo_define:nnnnnN
	{ \char_generate:nn { #1 } { 11 } }
	{ \char_generate:nn { #1 } { 11 } }
	{ scr }
	{ tl }
	{ \scrtl }
	\use:n
	
	\bamboo_define:nnnnnN
	{ \char_generate:nn { #1 } { 11 } }
	{ \char_generate:nn { #1 } { 11 } }
	{ }
	{ dt }
	{ \dot}
	\use:n
}
\int_step_inline:nnn { `a } { `z }
{
	\bamboo_define:nnnnnN
	{ \char_generate:nn { #1 } { 11 } }
	{ \char_generate:nn { #1 } { 11 } }
	{ bl }
	{ }
	{ \mymathbold }
	\use:n
	
	\bamboo_define:nnnnnN
	{ \char_generate:nn { #1 } { 11 } }
	{ \char_generate:nn { #1 } { 11 } }
	{ }
	{ hat }
	{ \hat }
	\use:n
	
	\bamboo_define:nnnnnN
	{ \char_generate:nn { #1 } { 11 } }
	{ \char_generate:nn { #1 } { 11 } }
	{ }
	{ br }
	{ \bar }
	\use:n
	
	\bamboo_define:nnnnnN
	{ \char_generate:nn { #1 } { 11 } }
	{ \char_generate:nn { #1 } { 11 } }
	{ }
	{ tl }
	{ \tilde }
	\use:n                            
	
	\bamboo_define:nnnnnN
	{ \char_generate:nn { #1 } { 11 } }
	{ \char_generate:nn { #1 } { 11 } }
	{ }
	{ dt }
	{ \dot}
	\use:n
}
\clist_map_inline:nn
{
	Gamma,Delta,Theta,Lambda,Xi,Pi,Sigma,Phi,Psi,Omega
}
{
	\bamboo_define:nnnnnN
	{ #1 }
	{ #1 }
	{ bl }
	{ }
	{ \mymathbold }
	\use:c
	
	\bamboo_define:nnnnnN
	{ #1 }
	{ #1 }
	{ op }
	{ }
	{ \op }
	\use:c
	
	\bamboo_define:nnnnnN
	{ #1 }
	{ #1 }
	{ scr }
	{ }
	{ \mathcal }
	\use:c
	
	\bamboo_define:nnnnnN
	{ #1 }
	{ #1 }
	{ }
	{ hat }
	{ \widehat }
	\use:c
	
	\bamboo_define:nnnnnN
	{ #1 }
	{ #1 }
	{ }
	{ br }
	{ \overline }
	\use:c
	
	\bamboo_define:nnnnnN
	{ #1 }
	{ #1 }
	{ }
	{ tl }
	{ \widetilde }
	\use:c
	
	\bamboo_define:nnnnnN
	{ #1 }
	{ #1 }
	{ }
	{ dt }
	{ \dot}
	\use:c
	
}
\clist_map_inline:nn
{
	alpha,beta,gamma,delta,epsilon,zeta,eta,theta,iota,kappa,
	lambda,mu,nu,xi,pi,rho,sigma,tau,phi,chi,psi,omega
}
{
	\bamboo_define:nnnnnN
	{ #1 }
	{ #1 }
	{ bl }
	{ }
	{ \mymathbold }
	\use:c
	
	\bamboo_define:nnnnnN
	{ #1 }
	{ #1 }
	{ }
	{ hat }
	{ \hat }
	\use:c
	
	\bamboo_define:nnnnnN
	{ #1 }
	{ #1 }
	{ }
	{ br }
	{ \bar }
	\use:c
	
	\bamboo_define:nnnnnN
	{ #1 }
	{ #1 }
	{ }
	{ tl }
	{ \tilde }
	\use:c
	
	\bamboo_define:nnnnnN
	{ #1 }
	{ #1 }
	{ }
	{ dt }
	{ \dot}
	\use:c
}

\ExplSyntaxOff

\DeclareMathOperator{\E}{\mathbf{E}}
\renewcommand{\P}{\operatorname{\mathbf{P}}}

\newcommand{\tr}{\operatorname{tr}}
\newcommand{\argmin}{\operatornamewithlimits{arg~min}}
\newcommand{\argmax}{\operatornamewithlimits{arg~max}}
\newcommand{\diag}{\operatorname{diag}}

\DeclarePairedDelimiter{\norm}{\lVert}{\rVert}

\DeclarePairedDelimiter{\abs}{\lvert}{\rvert}

\DeclarePairedDelimiter{\braces}{\{}{\}}
\DeclarePairedDelimiter{\parens}{(}{)}
\DeclarePairedDelimiter{\brackets}{[}{]}
\DeclarePairedDelimiterX{\ip}[2]{\langle}{\rangle}{#1,#2}

\DeclarePairedDelimiterXPP{\normsub}[2]{}{\lVert}{\rVert}{_{#2}}{#1}
\DeclarePairedDelimiterXPP{\ipsub}[3]{}{\langle}{\rangle}{_{#3}}{#1,#2}

\DeclarePairedDelimiterXPP{\ipHS}[2]{}{\langle}{\rangle}{_{\mathrm{HS}}}{#1, #2}
\DeclarePairedDelimiterXPP{\normHS}[1]{}{\lVert}{\rVert}{_{\mathrm{HS}}}{#1}

\DeclarePairedDelimiterXPP{\ipF}[2]{}{\langle}{\rangle}{_{\mathrm{F}}}{#1, #2}
\DeclarePairedDelimiterXPP{\normF}[1]{}{\lVert}{\rVert}{_{\mathrm{F}}}{#1}

\DeclarePairedDelimiterXPP{\dkl}[2]{\operatorname{D_{KL}}}{(}{)}{}{#1 \: \delimsize\Vert \: #2}

\DeclarePairedDelimiterXPP{\restr}[2]{}{{}}{\vert}{_{#2}}{#1}

\newcommand{\ones}{\mathbf{1}}

\newcommand{\R}{\mathbf{R}}
\newcommand{\C}{\mathbf{C}}
\newcommand{\Z}{\mathbf{Z}}

\newcommand{\spn}{\operatorname{span}}

\newcommand{\normaldist}{\operatorname{\mathcal{N}}}

\newcommand{\negqquad}{\hspace{-2em}}

\usepackage{bookmark}
\usepackage{color}
\usepackage[capitalize]{cleveref}

\usepackage{hyperref}

\DeclarePairedDelimiterXPP{\opnorm}[1]{}{\lVert}{\rVert}{_{\ell_2}}{#1}
\DeclarePairedDelimiterXPP{\nucnorm}[1]{}{\lVert}{\rVert}{_{*}}{#1}
\DeclarePairedDelimiterXPP{\normD}[1]{}{\lVert}{\rVert}{_D}{#1}
\DeclarePairedDelimiterXPP{\ipD}[2]{}{\langle}{\rangle}{_{D}}{#1, #2}

\newcommand\ddiag{\operatorname{ddiag}}

\newcommand{\SOgrp}{\mathrm{SO}}
\newcommand{\st}{\text{ s.t.\ }}
\newcommand{\ER}{Erd\H{o}s--R\'enyi}
\newcommand{\dmin}{d_{\mathrm{min}}}

\newcommand{\tansp}{T_Y}
\newcommand{\Ptansp}{\scrP_{\tansp}}
\newcommand{\Ponepmat}{P_{\ones}^\perp}
\newcommand{\PDonepmat}{P_{d}^\perp}
\newcommand{\real}{\operatorname{Re}}

\newcommand{\sphdim}{r}
\newcommand{\Sph}{\scrS}

\newcommand{\Sdm}{\Sph^{\sphdim-1}}

\newcommand{\rew}{Rakoto Endor and Waldspurger}

\newcommand{\longtitle}{Benign landscapes for synchronization on spheres\\ via normalized Laplacian matrices}
\newcommand{\shorttitle}{Benign landscapes for synchronization on spheres}
\newcommand{\fundingack}{This work was supported by the Swiss State Secretariat for Education, Research, and Innovation (SERI) under contract number MB22.00027 during the author's time at EPFL.}
\newcommand{\myaddress}{CERMICS, ENPC, Institut Polytechnique de Paris,	CNRS, Marne-la-Vallée, France}

\ifSIAM
\title{\longtitle%
	\thanks{Submitted to the editors March 24, 2025.%
	\funding{\fundingack{}}}}

\author{Andrew D.\ McRae%
	\thanks{\myaddress{} (\email{andrew.mcrae@enpc.fr}).}%
}

\headers{\shorttitle}{A. D. McRae}

\else
\hypersetup{
	bookmarks=true,
	colorlinks,
	linkcolor=blue,
	citecolor=blue,
	urlcolor=blue,
}
\newtheorem{theorem}{Theorem}[section]
\newtheorem{corollary}[theorem]{Corollary}

\title{\longtitle}

\author{Andrew D.\ McRae\thanks{\myaddress{} (\texttt{andrew.mcrae@epfl.ch}). \fundingack{}}}
\fi

\begin{document}

\maketitle

\begin{abstract}
	We study the nonconvex optimization landscapes of synchronization problems on spheres. First, we present new results for the statistical problem of synchronization over the two-element group $\mathbf{Z}_2$. We consider the nonconvex least-squares problem with $\mathbf{Z}_2 = \{\pm 1\}$ relaxed to the unit sphere in $\mathbf{R}^r$ for $r \geq 2$; for several popular models,	including graph clustering under the binary stochastic block model, we show that, for any $r \geq 2$, every second-order critical point recovers the ground truth in the asymptotic regimes where exact recovery is information-theoretically possible. Such statistical optimality via spherical relaxations had previously only been shown for (potentially arbitrarily) larger relaxation dimension $r$. Second, we consider the global synchronization of networks of coupled oscillators under the (homogeneous) Kuramoto model. We prove new and optimal asymptotic results for random signed networks on an Erd\H{o}s--R\'enyi graph, and we give new and simple proofs for several existing state-of-the-art results. Our key tool is a deterministic landscape condition that extends a recent result of Rakoto Endor and Waldspurger. This result says that, if a certain problem-dependent Laplacian matrix has small enough condition number, the nonconvex landscape is benign. Our extension allows the condition number to include an arbitrary diagonal preconditioner, which gives tighter results for many problems. We show that, for the synchronization of Kuramoto oscillator networks on nearest-neighbor circulant graphs as studied by Wiley, Strogatz, and Girvan, this condition is optimal. We also prove a natural complex extension that may be of interest for synchronization on the special orthogonal group $\operatorname{SO}(2)$.
\end{abstract}

\section{Introduction}
This paper considers optimization problems of the form
\begin{equation}
	\label{eq:opt_general}
	\max_{y_1, \dots, y_n \in \Sdm}~\sum_{i,j=1}^n C_{ij} \ip{y_i}{y_j}
	= \max_{Y \in \R^{n \times \sphdim}}~\ip{C}{Y Y^T} \st \diag(Y Y^T) = \ones,
\end{equation}
where $\Sdm$ denotes the unit sphere in $\R^\sphdim$, $C = [C_{ij}]$ is a symmetric $n \times n$ matrix,
$\ip{\cdot}{\cdot}$ denotes the standard Euclidean inner product (applied elementwise if the inputs are matrices, in which case it is often known as the trace or Frobenius inner product), $\diag(\cdot)$ applied to a square matrix extracts the diagonal as a vector, and $\ones$ denotes the all-ones vector in $\R^n$.

For $\sphdim \geq 2$, the problem \eqref{eq:opt_general} is smooth but nonconvex, so, in general, it can have spurious (non-global) local optima.
However, for certain choices of cost matrix $C$ arising in applications,
this nonconvex problem can have a benign landscape: this means that every local optimum (or, even stronger, every \emph{second-order critical point} where, in a Riemannian geometry sense, the gradient is zero and the Hessian is negative semidefinite) is, in fact, a global optimum.

We will specifically consider instances of \eqref{eq:opt_general} arising in two problems with quite different motivation but remarkably similar mathematical properties.

\subsection{$\Z_2$ synchronization and graph clustering}
The first application is statistical in nature.
Suppose, given a vector $z \in \{ \pm 1 \}^n$ that we want to estimate,
we make (potentially noisy) measurements of the relative signs between some pairs of elements in $z$:
we observe $R_{ij} \approx z_i z_j$ for certain (unordered) pairs $(i,j)$.
The problem of estimating $z$ from such measurements is called $\Z_2$ synchronization (due to the identification of $\{\pm 1\}$ under multiplication with the two-element group $\Z_2$) and arises in applications such as graph clustering (where the elements of $z$ represent graph vertex labels) and sensor network localization (where the signs represent reflections in space),
and it is one of the simplest instances of the broader problem of \emph{group synchronization}.
See \cite{Abbe2014,Cucuringu2012a,Bansal2004,McRae2025} for further background and references.
A least-squares estimate of $z$ is
\[
	\zhat \in \argmin_{x \in \{ \pm 1 \}^n} \sum_{(i,j) \in E} (R_{ij} - x_i x_j)^2 = \argmax_{x \in \{ \pm 1 \}^n} \sum_{(i,j) \in E} R_{ij} x_i x_j,
\]
where $E$ is the set of unordered pairs $(i,j)$ for which we have measurements.
This is precisely in the form of \eqref{eq:opt_general} in the case $\sphdim = 1$ (as $\Sph^0 = \{ \pm 1 \}$), where $C$ is constructed from the available measurements $R_{ij}, (i,j) \in E$.
Larger choices of $\sphdim$ correspond to \emph{relaxing} the problem,
that is, adding variables in the hope of making it easier to solve.
In particular, choosing $\sphdim \geq 2$ allows the use of continuous optimization algorithms such as gradient descent and trust-region methods.
The choice $\sphdim = 2$ was introduced for $\Z_2$ synchronization by Bandeira et~al.\ \cite{Bandeira2016a};
the same type of relaxation was introduced for the max-cut problem (which also has the form of \eqref{eq:opt_general} for $\sphdim = 1$) by Burer et~al.\ \cite{Burer2002}.
Relaxations with larger $\sphdim$ have been studied more recently by Ling \cite{Ling2025}, McRae et~al.\ \cite{McRae2025}, and \rew{} \cite{Endor2026}.
For several popular (random) models of this problem, 
the optimal requirements (in terms of number of measurements and amount of noise) for recovering $z$ exactly (in an asymptotic sense for large problem sizes) are well-established.
In particular, we consider the following models (see \Cref{sec:intro_appls} for more details and background):
\begin{itemize}
	\item \textbf{Gaussian noise:} For every pair $(i,j)$, we observe $z_i z_j + \xi_{ij}$, where $\{\xi_{ij}\}_{i<j}$ are independent and identically distributed (i.i.d.) $\normaldist(0, \sigma^2)$ random variables for some $\sigma \geq 0$. The exact recovery threshold (asymptotically, as the problem size $n \to \infty$) is $\sigma = \sqrt{\frac{n}{2 \log n}}$.
	\item \textbf{\ER{} random graph with Bernoulli noise:} For each pair $(i,j)$, we measure $z_i z_j$ with probability $p > 0$. This gives us an \ER{} random \emph{measurement graph}.
	Furthermore, for each measurement we make, the sign is randomly flipped with probability $\frac{1 - \delta}{2}$ for some parameter $\delta \in [0,1]$ that represents the signal strength.
	The exact recovery threshold is $\frac{np}{\log n}( 1 - \sqrt{1 - \delta^2} ) = 1$.
	\item \textbf{Stochastic block model (SBM) in graph clustering:} In this model, the elements of $z$ represent community (cluster) labels, and we want to estimate these from a random graph whose $(i,j)$ edge is present with probability $p$ if $z_i = z_j$ and with probability $q \neq p$ if $z_i \neq z_j$.
	The exact cluster recovery threshold is $\frac{n}{\log n} (\sqrt{p} - \sqrt{q})^2 = 2$.
\end{itemize}
In particular, it is known (again, see \Cref{sec:intro_appls}) that the semidefinite relaxation (equivalent to $\sphdim = n$; see \Cref{sec:intro_opt}) exactly recovers $z$ up to these thresholds, but the number of optimization variables is then proportional to $n^2$, which may be computationally impractical for large $n$.
The recent results of \cite{McRae2025} (and, for certain problems, \cite{Endor2026}) show that we can obtain the same statistical performance with smaller values of $\sphdim$ (in the sense that every second-order critical point $Y$ of \eqref{eq:opt_general} for appropriate $C$ satisfies $Y Y^T = z z^T$),
but, the closer we are to the threshold, the larger $\sphdim$ needs to be for these results to apply.

In this work, we prove that, for the above problems, the smallest continuous relaxation ($\sphdim = 2$) is sufficient to obtain asymptotic exact recovery all the way up to the optimal threshold.
This shows that optimization over $\approx n$ variables suffices to recover $z$ exactly whenever it is statistically possible (in an asymptotic sense);
this is an improvment over the $\approx n^2$ variables required for the semidefinite relaxation and the potentially large multiple of $n$ required by the results of \cite{McRae2025}.

\subsection{(Signed) Kuramoto oscillator networks}
\label{sec:intro_kuramoto}
A second application comes from the study of dynamical systems.
A simple (homogeneous) version of the Kuramoto model \cite{Kuramoto1975} for a network of coupled oscillators is the following:
we have $n$ time-varying angles (oscillator phases) $\theta_1(t), \dots, \theta_n(t)$
that evolve simultaneously according to the system of ordinary differential equations
\begin{equation}
	\label{eq:kuramoto}
	\frac{d \theta_i}{d t} = K \sum_{j=1}^n A_{ij} \sin(\theta_j - \theta_i),
\end{equation}
where $A = [A_{ij}]$ is a symmetric $n \times n$ coupling matrix (often taken to be the adjacency matrix of a graph), and $K > 0$ is a global coupling constant ($K$ only affects the time scaling and not the stability or long-term convergence of the system, so it will play no role in the rest of this paper).
It is clear that a family of equilibria of \eqref{eq:kuramoto} is $\theta_1 = \cdots = \theta_n \mod 2\pi$.
Furthermore, when $A$ is the adjacency matrix of a (potentially weighted) connected graph,
one can easily check that these equilibria are \emph{stable}.\footnote{The Jacobian at such points is the negative graph Laplacian, which has strictly negative eigenvalues except in the trivial direction corresponding to shifting all angles equally.}
But are these the only stable equilibria?
Or, even stronger, is it the case that the dynamical system will converge to such a ``synchronized'' state from (almost) every starting point?
In general, the answer is no, but for certain coupling matrices $A$ the answer is indeed yes;
we say that the oscillator network associated to $A$ \emph{globally synchronizes}.%
\footnote{If $A$ is the adjacency matrix of a graph $G$, we may say that $G$ is globally synchronizing.}
The connection to the optimization problem \eqref{eq:opt_general} is that the dynamical system \eqref{eq:kuramoto} is, up to time scaling and with the angle parametrization of the unit circle, precisely the (maximizing) Riemannian gradient flow of \eqref{eq:opt_general} with $\sphdim = 2$ and $C = A$.
As observed by Ling et~al.\ \cite{Ling2019}, to guarantee global synchronization, it suffices\footnote{To make this fully rigorous, we also need the fact that the objective function and constraint set are analytic to rule out certain pathological behavior. See, for example, \cite{Markdahl2020} or \cite{Geshkovski2025} for such an argument.} to show that every second-order critical point $Y$ of \eqref{eq:opt_general} represents a synchronized state in the sense that $Y Y^T = \ones \ones^T$ (in the angle formulation, this is equivalent to $\cos(\theta_i - \theta_j) = 1$ for all $i,j$).
Determining which graphs globally synchronize is a highly active area of research;
recent results consider random graphs, dense graphs, expander graphs, and various structured graphs \cite{Abdalla2026,Abdalla2024,Kassabov2021,DeVita2025,Sclosa2023,Yoneda2021}.
Nevertheless, many issues and open questions remain; in particular, we seek in this paper to address the following:
\begin{itemize}
	\item Most state-of-the-art results only consider the case where $A$ is the adjacency matrix of a graph, that is, the coupling coefficients are nonnegative.
	What happens when a few coefficients are negative (i.e., certain oscillators repulse each other)?
	Intuitively (and empirically), a few ``negative edges'' should not prevent synchronization if there are enough positive edges to compensate, but many of the current proof techniques cannot handle any negative coefficients.
	\item Relatedly, many recent results use specialized analysis techniques which, despite some similarities, do not transfer in a straightforward way from one class of network to another.
	Can we unify these results under a common approach?
\end{itemize}
We present a simple (and, for certain problem types, optimal) condition on the coupling matrix $A$ guaranteeing that the oscillator network with dynamics \eqref{eq:kuramoto} is globally synchronizing.
From this, we recover many recent state-of-the-art results with minimal additional effort.
In addition, our synchronization condition allows for the presence of negative (repulsive) coupling,
and we obtain new (and optimal) results for synchronization on an \ER{} graph with some (randomly) negative edges.

\subsection{Paper outline and additional notation}
The rest of this paper is organized as follows:
\begin{itemize}
	\item \Cref{sec:intro_opt} describes optimality conditions for \eqref{eq:opt_general} via its (convex) semidefinite relaxation
	and gives a general result (\Cref{thm:sync_cn}, which is an extension of the main result of \rew{} \cite{Endor2026}) showing that, if a certain Laplacian matrix (normalized by an arbitrary diagonal preconditioner) has a smaller condition number than the sphere dimension $\sphdim$, the problem \eqref{eq:opt_general} has a benign landscape.
	We also provide an extension (\Cref{thm:sync_cn_cplx}) to the complex case, which may be of interest for other applications.
	
	\item \Cref{sec:intro_appls} describes the $\Z_2$ synchronization and Kuramoto oscillator network applications in more detail and states our main results for these applications.
	In particular, \Cref{sec:intro_gaussian,sec:intro_bern,sec:intro_sbm} describe in detail the three $\Z_2$ synchronization models mentioned above,
	and their main results, \Cref{thm:z2_gauss,thm:z2_bern,thm:sbm}, show that the nonconvex optimization approach of \eqref{eq:opt_general} with any $r \geq 2$ asymptotically achieves exact recovery all the way up to the optimal threshold.
	\Cref{thm:z2_bern} can also be interpreted in terms of signed Kuramoto oscillator networks and gives a new and optimal result.
	\Cref{sec:intro_kuramoto_further} discusses several additional popular Kuramoto oscillator network models and shows that \Cref{thm:sync_cn}, restated and shown to be optimal for Kuramoto oscillator networks as \Cref{cor:kuramoto_det},
	provides a simple proof for many existing state-of-the-art results.
	
	\item \Cref{sec:LD_condnum,sec:proofs_appls,sec:proof_condnum} contain the main technical proofs of this paper.
	In particular, \Cref{sec:LD_condnum} gives a bound (\Cref{thm:normalized_cn}) on the condition number of a particular normalized Laplacian matrix; we use this extensively in our applications, but it may also be of independent interest.
	\Cref{sec:proofs_appls} gives problem-specific proofs for the $\Z_2$ synchronization models and for one of the Kuramoto oscillator network models.
	\Cref{sec:proof_condnum} gives complete proofs of the deterministic landscape results \Cref{thm:sync_cn,thm:sync_cn_cplx}.
\end{itemize}

We will use the following notation throughout this paper.
If $A$ and $B$ are matrices, $A \circ B$ denotes their elementwise (Hadamard) product, and $A^{\circ 2}$ denotes the elementwise (Hadamard) square of $A$.
We will denote the (real) eigenvalues of a Hermitian $n \times n$ matrix $A$ by $\lambda_1(A) \leq \cdots \leq \lambda_n(A)$ in increasing order.
If $v$ is a vector, $\norm{v}$ and $\norm{v}_1$ respectively denote its Euclidean ($\ell_2$) and $\ell_1$ norms.
If $A$ is a matrix, $\opnorm{A}$ and $\normF{A}$ respectively denote its $\ell_2$ operator norm and Frobenius (elementwise Euclidean) norm.
If $A$ is any real (resp.\ complex) $n \times n$ matrix, $\diag(A) \in \R^n$ (resp.\ $\C^n$) denotes the diagonal of $A$ extracted as a vector,
whereas, if $x \in \R^n$ (or $\C^n$), $\diag(x)$ denotes the diagonal $n \times n$ matrix with the elements of $x$ on the diagonal.
$\ddiag(A)$ is the diagonal matrix with the same diagonal as $A$ but with the remaining elements set to zero.
$\tr(A)$ is the trace (sum of the diagonal elements) of $A$.
$\real(x)$ denotes the real part of complex $x$, applied elementwise if $x$ is a vector or matrix.
We write $a \lesssim b$ (equivalently, $b \gtrsim a$) if we have $a \leq C b$ for some unspecified constant $C > 0$.
We write $a \approx b$ if $a \lesssim b$ and $a \gtrsim b$.
We also use the big-oh notation $O(b)$ to denote an unspecified quantity $a$ with $\abs{a} \lesssim b$.
We denote $a \wedge b = \min\{ a, b \}$.

\section{Optimality conditions and deterministic landscape results}
\label{sec:intro_opt}
To understand the optimality of (candidate) solutions to \eqref{eq:opt_general},
it is useful to consider the convex \emph{semidefinite relaxation},
as optimality conditions for convex problems are well understood.
Noting that the variable $Y$ in \eqref{eq:opt_general} only appears in the objective function and constraint via $Y Y^T$, which is a positive semidefinite (PSD) matrix of rank at most $\sphdim$,
we can relax the rank constraint\footnote{Setting $\sphdim = n$ in \eqref{eq:opt_general} makes the two problems equivalent. The problem \eqref{eq:opt_general} is still nonconvex, but it turns out that, when $\sphdim > \frac{n}{2}$, the landscape is benign for any $C$ \cite[Cor.~5.11]{Boumal2019}: for every second-order critical point $Y$, $X = Y Y^T$ is a global optimum of \eqref{eq:opt_sdp}.}  to obtain the convex semidefinite program (SDP)
\begin{equation}
	\label{eq:opt_sdp}
	\max_{X \succeq 0}~\ip{C}{X} \st \diag(X) = \ones.
\end{equation}
This is often called the ``max-cut SDP,'' as a semidefinite relaxation of this form was first introduced for the max-cut problem \cite{Goemans1995}.
In fact, reversing the relaxation process, the nonconvex problem \eqref{eq:opt_general} can be seen as a rank-$\sphdim$ Burer-Monteiro factorization \cite{Burer2003} of \eqref{eq:opt_sdp}.

How might we show that a candidate solution $Y$ to \eqref{eq:opt_general} is globally optimal?
As relaxing the problem can only increase the optimal value, a clearly \emph{sufficient} condition is for $X = Y Y^T$ to be an optimal solution to the semidefinite relaxation \eqref{eq:opt_sdp}.

How, in turn, might we show that a feasible point $X$ of \eqref{eq:opt_sdp} is optimal?
Standard convex duality theory says that $X$ is optimal if there exists a diagonal matrix $\Lambda$ such that the matrix $S \coloneqq \Lambda - C$ satisfies $S \succeq 0$ and $S X = 0$.
Indeed, if $X'$ is another feasible point, we must have
\begin{align*}
	\ip{C}{X} - \ip{C}{X'}
	&= \ip{S}{X'} - \underbrace{\ip{S}{X}}_{= 0} + \underbrace{\ip{\Lambda}{X - X'}}_{= 0} \\
	&= \ip{S}{X'} \\
	&\geq 0.
\end{align*}
We have $\ip{\Lambda}{X - X'} = 0$ because $\Lambda$ is diagonal and $X$ and $X'$ are equal on their diagonals.
The last inequality is due to the fact that $S \succeq 0$ (by assumption) and $X' \succeq 0$ (by feasibility).

We now use the fact that we hope for a particular rank-1 solution.%
\footnote{This optimality analysis hinges on knowing the global optimum beforehand. The applications we consider allow for exact recovery/synchronization, so for analysis purposes we can simply plug in the ground truth in the hope of showing it is, in fact, the global optimum. Generalizations of this problem, such as higher-dimensional orthogonal group synchronization, do not allow exact recovery when there is measurement noise, which makes landscape analysis more difficult---see, for example, \cite{McRae2024,Ling2025}.}
Recalling the applications described above, we want to show that, for a fixed $z \in \{ \pm 1 \}^n$, $X = z z^T$ is a solution to \eqref{eq:opt_sdp}.
In this case, the optimality condition $S X = 0$ is equivalent to $S z = 0$ which, in turn, implies
\[
	\Lambda = \diag(z) \diag(C z) = \ddiag(C z z^T) \quad \Longrightarrow \quad S = L(z) \coloneqq \ddiag(C z z^T) - C.
\]
By construction, $L(z) z = 0$; it is a valid optimality certificate for $X = z z^T$ if, in addition, $L(z) \succeq 0$.
We call $L(z)$ a Laplacian matrix as, in the case that $C$ is a graph adjacency matrix and $z = \ones$,
it is the ordinary graph Laplacian.

Because $L(z) z = 0$, $L(z)$ has at least one zero eigenvalue.
Then, if $L(z) \succeq 0$, we must have $\lambda_1(L(z)) = 0$.
The next eigenvalue $\lambda_2(L(z))$ (known as the spectral gap or, in the graph Laplacian case, the algebraic connectivity or the Fiedler value) measures how robust the solution $X = z z^T$ is;
$\lambda_2(L(z)) > 0$ ensures that this solution is \emph{unique}.

Therefore, in the case that $\lambda_2(L(z)) > 0$,
\emph{every} global solution $Y$ to the nonconvex problem \eqref{eq:opt_general} must satisfy $Y Y^T = z z^T$ (or, equivalently, $Y = z v^T$ for some unit-norm $v \in \R^r$).
This does not, however, rule out spurious local optima.

The remarkable recent results of Ling \cite{Ling2025} and \rew{} \cite{Endor2026} tell us that we can know something about the nonconvex landscape of \eqref{eq:opt_general} from the \emph{condition number} of $L(z)$.
In particular, the more recent and stronger result of \cite{Endor2026} says that, if the sphere dimension satisfies
\[
	\sphdim > \frac{\lambda_n(L(z))}{\lambda_2(L(z))},
\]
then every second-order critical point $Y$ of \eqref{eq:opt_general} is a global optimum satisfying $Y Y^T = z z^T$.
Conversely, they show that for any $\sphdim \geq 2$, there is a problem instance with a spurious second-order critical point and $\frac{\lambda_n(L(z))}{\lambda_2(L(z))} = \sphdim$.
We provide a somewhat different class of counterexamples in the case $\sphdim = 2$ in \Cref{sec:intro_kuramoto_further}.

The following theorem is a generalization of the result of \cite{Endor2026}:
it applies with any \emph{diagonal preconditioner}, which, as we will see for certain applications, can considerably improve the condition number.
\begin{theorem}
	\label{thm:sync_cn}
	Let $z \in \{ \pm 1 \}^n$, and suppose $L(z) \coloneqq \ddiag(C z z^T) - C \succeq 0$ with $\lambda_2(L(z)) > 0$.
	Let $D$ be any (strictly) positive definite diagonal matrix, and set $L_D(z) \coloneqq D^{-1/2} L(z) D^{-1/2}$.
	Then, if
	\[
		\sphdim > \frac{\lambda_n(L_D(z))}{\lambda_2(L_D(z))},
	\]
	every second-order critical point $Y$ of \eqref{eq:opt_general} is a global optimum satisfying $Y Y^T = z z^T$.
\end{theorem}
For completeness, we provide a full proof in \Cref{sec:proof_condnum_real}.
Our proof is inspired by but significantly shorter than that of \rew{} \cite{Endor2026}.

The motivation for the diagonal preconditioner is that, for many of the applications we consider, the condition number of $L(z)$ depends on the largest and smallest diagonal elements of $\ddiag(C z z^T)$.
Normalizing by a diagonal preconditioner (in this paper, we use $\ddiag(C z z^T)$ itself) reduces this issue.

Apart from \cite{Endor2026}, the closest results to which we can compare \Cref{thm:sync_cn} are \cite[Thm.~2.1]{Ling2025}, which requires
\[
	\sphdim \geq 2\frac{\lambda_n(L(z))}{\lambda_2(L(z))} + 1,
\]
and \cite[Thm.~1]{McRae2025}, which is more specialized to a particular structure of $C$ but also requires $\sphdim \geq 3$.
\Cref{thm:sync_cn} allows $\sphdim$ as small as $2$.

On the other hand, a useful feature of the results of \cite{Ling2025,McRae2025} is their robustness to \emph{monotone adversaries}.
This means that adding a matrix $\Delta$ to the cost matrix $C$, where the entries of $\Delta$ satisfy $\Delta_{ij} z_i z_j \geq 0$,
does not negatively affect the theoretical guarantees in those results, just as it does not change the global optimum of \eqref{eq:opt_general}.
See \cite[Sec.~5]{McRae2025} for more details.
\Cref{thm:sync_cn} is \emph{not} robust to monotone adversaries in the case $\sphdim = 2$ (one can verify that, similarly to \cite{Ling2025,McRae2025}, our result is robust when $\sphdim \geq 3$).
For example, given any cost matrix $C$ to which \Cref{thm:sync_cn} applies with $\sphdim = 2$ and $z = \ones$,
consider the problem with cost matrix $C' = C + \alpha A$,
where $\alpha \geq 0$ is some scale parameter, and $A$ is the adjacency matrix of a cycle graph.
For sufficiently large $\alpha$, the landscape becomes arbitrarily close to that of synchronization on a cycle graph,
which is known \cite{Wiley2006} to have spurious local optima.
This resembles the observation of Ling et~al.\ \cite{Ling2019} that adding edges to globally synchronizing Kuramoto oscillator networks can give rise to spurious stable states.

\subsection*{Complex extension}
\Cref{thm:sync_cn} and its proof generalize in a straightforward way to the complex case.
We state and prove this result but do not, in the present paper, explore its implications.
One interesting potential application is \emph{synchronization of rotations}, as the complex unit circle is the two-dimensional special orthogonal group (rotation group) $\SOgrp(2)$. This problem has many applications in areas such as robotics and computer vision; see, for example, the recent papers \cite{Dellaert2020,Doherty2022} for more background and references.
The proof of \Cref{thm:sync_cn} (here and in \cite{Endor2026}) does not generalize in an obvious way to higher-dimensional orthogonal groups as does that of, for example, \cite{Ling2025},
so it is interesting that it nevertheless does apply (via this complex extension) to the specific group $\SOgrp(2)$.

We now let $C$ be an $n \times n$ complex Hermitian matrix, and we consider the optimization problem
\begin{equation}
	\label{eq:opt_cplx}
	\max_{Y \in \C^{n \times \sphdim}}~\ip{C}{Y Y^*} \st \diag(Y Y^*) = \ones.
\end{equation}
Criticality conditions for the complex variable $Y = U + iV$ should be considered in terms of its real and imaginary parts, that is, the variable $(U, V)$.

To show that the semidefinite relaxation has rank-1 optimal solution $X = z z^*$ for some $z \in \C^n$ (with, necessarily, $\abs{z_i} = 1$ for all $i$), the dual certificate matrix is now
\[
	L(z) \coloneqq \real(\ddiag(C z z^*)) - C.
\]
Again, we need $L(z) z = 0$
and $L(z) \succeq 0$ to certify optimality.
Note that, unlike in the real case, $L(z)z = 0$ does not automatically hold by construction (in fact, it requires $\diag(C z z^*)$ to be real).

With this in place, we can state the complex extension of \Cref{thm:sync_cn} (proved in \Cref{sec:proof_condnum_cplx}):
\begin{theorem}
	\label{thm:sync_cn_cplx}
	Let $z \in \C^n$ with $\abs{z_i} = 1$ for all $i = 1, \dots, n$, and suppose $L(z) \coloneqq  \real(\ddiag(C z z^*)) - C$ satisfies $L(z) z = 0$ and $L(z) \succeq 0$ with $\lambda_2(L(z)) > 0$.
	Let $D$ be any (strictly) positive definite diagonal matrix, and set $L_D(z) \coloneqq D^{-1/2} L(z) D^{-1/2}$.
	Then, if
	\[
		2 \sphdim > \frac{\lambda_n(L_D(z))}{\lambda_2(L_D(z))},
	\]
	every second-order critical point $Y$ of \eqref{eq:opt_cplx} is a global optimum satisfying $Y Y^* = z z^*$.
\end{theorem}
Note that if $C$ and $z$ are, in fact, real, this directly reduces to \Cref{thm:sync_cn} under the identification of the unit sphere in $\C^\sphdim$ with $\Sph^{2\sphdim-1}$.

\section{Applications}
\label{sec:intro_appls}
\subsection{$\Z_2$ synchronization with Gaussian noise}
\label{sec:intro_gaussian}
One of the simplest instances of the $\Z_2$ synchronization problem is where, given a ground truth $z \in \{\pm 1\}^n$,
we observe, for each (unordered) pair $(i,j)$, the noisy relative measurement $z_i z_j + \xi_{ij}$,
where $\xi_{ij} \sim \normaldist(0, \sigma^2)$ for some $\sigma \geq 0$,
and the noise is independent across all pairs $(i,j)$.
In matrix form, this gives
\begin{equation}
	\label{eq:z2_model_gauss}
	C = z z^T + \sigma W,
\end{equation}
where $W$ is a symmetric matrix of independent (up to symmetry) and identically distributed standard normal random variables (modulo the irrelevant diagonal elements and scaling, $W$ is a draw from the Gaussian orthogonal ensemble).
This model has been studied as a simplification of (and thus a stepping stone toward understanding) more complex models \cite{Bandeira2017,Javanmard2016}.

Under this model, we have the following asymptotic result, proved in \Cref{sec:proofs_z2}:
\begin{theorem}
	\label{thm:z2_gauss}
	Consider the model \eqref{eq:z2_model_gauss}.
	If the noise parameter $\sigma$ (which can depend on the problem size $n$) satisfies, for some $\epsilon > 0$ independent of $n$,
	\[
		\sigma \leq \sqrt{\frac{n}{(2 + \epsilon) \log n}},
	\]
	then, with probability $\to 1$ as $n \to \infty$,
	for all $\sphdim \geq 2$, every second-order critical point $Y$ of \eqref{eq:opt_general} with cost matrix $C$ satisfies $Y Y^T = z z^T$.
\end{theorem}
As shown by Bandeira \cite{Bandeira2018}, the threshold $\sqrt{\frac{n}{2 \log n}}$ is optimal for it to be possible to recover $z$ exactly by the maximum likelihood estimator,
and the semidefinite relaxation reaches this threshold.
However, previous results for the nonconvex approach \cite{McRae2025,Endor2026} required $\sphdim$ to be large to approach this threshold;
for $\sphdim = 2$, the best previous result \cite[Corollary 3.1]{Endor2026} is off by a factor of $3$ in terms of the bound on $\sigma$.
Bandeira et~al.\ \cite{Bandeira2016a} previously showed exact recovery with $r = 2$ under the stronger condition $\sigma \lesssim n^{1/6}$.

\subsection{$\Z_2$ synchronization/``signed Kuramoto'' on an \ER{} random graph with Bernoulli noise}
\label{sec:intro_bern}
We now consider another statistical model for $\Z_2$ synchronization with an \ER{} measurement graph and random sign flips.
This is motivated by the problem of clustering a graph with signed edges,
historically called \emph{correlation clustering} \cite{Bansal2004}.
The specific model we consider here was introduced by Abbé et~al.\ \cite{Abbe2014}, who called it the \emph{censored block model}.

For each (unordered) pair $(i,j)$, we observe the relative sign $z_i z_j$ with probability $p$,
but this measurement can be randomly flipped.
More precisely, the model for parameters $p, \delta \in [0, 1]$, is that, for each pair $(i,j)$ with $i < j$,
\begin{equation}
	\label{eq:z2_model_bern}
	C_{ij} = \begin{cases}
		z_i z_j & \text{ with probability } \frac{1 + \delta}{2} p, \\
		-z_i z_j & \text{ with probability } \frac{1 - \delta}{2} p, \text{ and} \\
		0 & \text{ with probability } 1 - p
	\end{cases}
\end{equation}
independently of the other pairs, and $C_{ji} = C_{ij}$ (again, the diagonal elements are irrelevant).
The parameter $p$ is the edge probability of the measurement graph,
while $\delta$ is a measure of the signal strength: if $\delta = 0$, the observed signs are completely random and uninformative, while if $\delta = 1$ there are no sign errors.

Alternatively, with $z = \ones$, the same random model for $C$ also has an interpretation for Kuramoto oscillator networks.
The coupling matrix $A = C$ corresponds to the adjacency matrix of an \ER{} random graph,
but some edges are randomly negative (i.e., the coupling is repulsive for those edges).
We will call this the \emph{random signed Kuramoto} model.
The case $p = 1$ was suggested by Ling \cite{Ling2025} for ``oscillators'' on higher-dimensional spheres.
Bandeira \cite[Conj.~9]{Bandeira2024oberwolfach} discusses the standard angular version and makes a conjecture which we resolve positively.

For this model, we have the following asmyptotic result, proved in \Cref{sec:proofs_z2}:
\begin{theorem}
	\label{thm:z2_bern}
	Consider the model \eqref{eq:z2_model_bern} with parameters $(p, \delta)$ that can  depend on the problem size $n$. Suppose, for some $\epsilon \in (0, 1)$ independent of $n$, that
	\[
		\frac{n p}{\log n} \parens*{ 1 - \sqrt{1 - (1 - \epsilon) \delta^2} } \geq 1.
	\]
	Then, with probability $\to 1$ as $n \to \infty$, the following hold:
	\begin{itemize}
		\item For any $\sphdim \geq 2$, every second-order critical point $Y$ of \eqref{eq:opt_general} with cost matrix $C$ satisfies $Y Y^T = z z^T$.
		\item Equivalently (for $\sphdim = 2$), the random signed Kuramoto oscillator network is globally synchronizing.
	\end{itemize}	
\end{theorem}
Note that in the case $\delta = 1$ (no random sign flips),
the condition reduces to $p \geq (1 + \epsilon) \frac{\log n}{n}$ (possibly for different $\epsilon > 0$),
which is the standard connectivity threshold for an \ER{} graph.
This recovers one of the main results of Abdalla et~al.\ \cite{Abdalla2026} (see \Cref{sec:kuramoto_er} for further related results and discussion).

Note also that
\[
	1 - \sqrt{1 - (1 - \epsilon) \delta^2} \geq \frac{1 - \epsilon}{2} \delta^2
\]
(for small $\delta$, this is a good approximation),
so we can strengthen the condition of \Cref{thm:z2_bern} to
\[
	\delta \geq \sqrt{\frac{2 (1 + \epsilon) \log n}{n p}}
\]
(again, $\epsilon > 0$ may be different).
In the case $p = 1$, this resolves the conjecture of Bandeira \cite[Conj.~9]{Bandeira2024oberwolfach}.

This is the optimal threshold for exact recovery of $z$ (or, in the language of the Kuramoto oscillator, for the network to synchronize even locally) \cite{Hajek2016a,Jog2015}.
As with \Cref{thm:z2_gauss}, it was known that the semidefinite relaxation reached this threshold \cite{Hajek2016a} and that, for large $\sphdim$, the nonconvex method approached this threshold (see \cite{McRae2025} for general $p$ and \cite{Endor2026} when $p = 1$),
but our result shows that we can reach this threshold even with $\sphdim = 2$.
For $\sphdim = 2$, the best previous result \cite[Cor.~3.2]{Endor2026} is off by a constant factor and only applies to the case $p = 1$.

\subsection{Stochastic block model}
\newcommand{\sbmcenterconst}{c}
\label{sec:intro_sbm}
With some additional effort, our results can be adapted to the (binary, symmetric) stochastic block model (SBM) in graph clustering.
This is a simple but theoretically important statistical model in the larger field of community detection.
See the surveys \cite{Fortunato2016,Moore2017,Abbe2018} for additional background. 

Let $z \in \{ \pm 1 \}^n$,
and assume, for simplicity, that $n$ is even and $\ip{z}{\ones} = 0$.
We observe an ordinary graph with adjacency matrix $A$,
where we assume that the presence of each edge is independent of the others and has probability (for $i < j$)
\begin{equation}
	\label{eq:sbm_model}
	\P(A_{ij} = 1) = \begin{cases}
		p & \text{if } z_i = z_j, \text{ and} \\
		q &\text{if } z_i \neq z_j,
	\end{cases}
\end{equation}
where $p \geq q \geq 0$.%
\footnote{The case $q > p$ can also be handled (it would require reversing the sign in our optimization problem), but for simplicity and brevity we omit this case.}
We want to recover $z$ from $A$.
We can easily calculate
\[
	\E A = \frac{p-q}{2} z z^T + \frac{p+q}{2} \ones \ones^T - p I_n.
\]
Attempting to solve \eqref{eq:opt_general} with $C = A$ directly would give an all-ones solution,
which is not interesting.
However, we can get around this by instead using the centered matrix $C = A - \sbmcenterconst \ones \ones^T$,
where $\sbmcenterconst$ is $\frac{p+q}{2}$ (if this is known) or some estimate of it (e.g., $c = \frac{1}{n^2} \ip{A}{\ones \ones^T}$).

For this model, we have the following result, proved in \Cref{sec:proofs_z2}.
\begin{theorem}
	\label{thm:sbm}
	Consider the model \eqref{eq:sbm_model} where the model parameters $p \geq q$ (possibly depending on the problem size $n$) satisfy, for some $\epsilon \in (0, 1/2)$ independent of $n$,
	\[
	\frac{n}{\log n} \parens*{ \sqrt{(1 - \epsilon)p + \epsilon q} - \sqrt{(1 - \epsilon)q + \epsilon p} }^2 \geq 2.
	\]
	Then, with probability $\to 1$ as $n \to \infty$, for $\sphdim \geq 2$, every second-order critical point of \eqref{eq:opt_general} with cost matrix $C = A - \sbmcenterconst \ones \ones^T$
	satisfies $Y Y^T = z z^T$,
	where either $\sbmcenterconst = \frac{p + q}{2}$ or $\sbmcenterconst = \frac{1}{n^2} \ip{A}{\ones \ones^T}$.
\end{theorem}

If we have $p = \frac{a \log n}{n}$, $q = \frac{b \log n}{n}$ for some $a \geq b \geq 0$,
we see that it suffices to have	$\sqrt{a} - \sqrt{b} > \sqrt{2}$.
Once again, this is the information-theoretically optimal threshold for exact recovery \cite{Abbe2016,Mossel2015}.
It was already shown that the semidefinite relaxation \cite{Hajek2016,Bandeira2018} and the nonconvex method for large $\sphdim$ \cite{McRae2025} reach this threshold,
but now we have shown that the nonconvex method reaches it even with $\sphdim = 2$.
Bandeira et~al.\ \cite{Bandeira2016a} earlier showed exact recovery with $\sphdim = 2$ under the stronger condition $(\sqrt{p} - \sqrt{q})^2 \gtrsim n^{-1/3}$.

\subsection{More on Kuramoto oscillator networks}
\label{sec:intro_kuramoto_further}
Restricted to two (standard) choices of preconditioner, \Cref{thm:sync_cn} implies a simple condition for Kuramoto oscillator networks to synchronize globally, which is the positive part of the following result:
\begin{theorem}
	\label{cor:kuramoto_det}
	Let $A$ be an $n \times n$ symmetric matrix.
	Assume the vector $A \ones$ has strictly positive entries,
	and set $D \coloneqq \diag(A \ones)$.
	Let $L = D - A$ and $\scrL = I_n - D^{-1/2} A D^{-1/2}$ be the ordinary and symmetric normalized Laplacian matrices.
	If $\lambda_2(L) > 0$ and either
	\[
		\frac{\lambda_n(L)}{\lambda_2(L)} < 2 \qquad \text{or} \qquad \frac{\lambda_n(\scrL)}{\lambda_2(\scrL)} < 2,
	\]
	then the Kuramoto oscillator network with dynamics \eqref{eq:kuramoto} with coupling matrix $A$ globally synchronizes.
	
	Conversely, for any $\epsilon > 0$,
	there is a connected, unweighted, regular graph $G$ with adjacency matrix $A$ such that
	\[
		\frac{\lambda_n(L)}{\lambda_2(L)} = \frac{\lambda_n(\scrL)}{\lambda_2(\scrL)} \leq 2 + \epsilon,
	\]
	and $G$ is not globally synchronizing.
\end{theorem}
The counterexamples proving the second part are provided by \Cref{thm:circulant} below.
In the rest of this section, we primarily focus on the case where $A$ is a graph adjacency matrix, and, therefore, $L$ and $\scrL$ are the ordinary and (symmetric) normalized graph Laplacians.
The result depending on $L$ is an immediate consequence of the result of \rew{} \cite{Endor2026};
the fact that this also works for the \emph{normalized} Laplacian $\scrL$ is, to our knowledge, new.

In the following sections, we consider how \Cref{cor:kuramoto_det} applies to several popular classes of graphs whose synchronization properties have been well studied. We do not obtain qualitatively new synchronization results (though \Cref{sec:expander} on regular expander graphs has some improved constants),
but we will see that, from \Cref{cor:kuramoto_det}, we can easily recover a number of known state-of-the-art results,
many of whose existing proofs are quite complex and specialized.

\subsubsection{Nearest-neighbor circulant graphs}
\label{sec:kuramoto_circ}
Before passing to more difficult examples,
we consider a simple class of graphs that has been well studied from the perspective of synchronization,
and we show that the condition of \Cref{cor:kuramoto_det} is optimal for these graphs.
These are the $k$-nearest-neighbor circulant graphs studied for synchronization by Wiley et~al.\ \cite{Wiley2006}.
For an integer $k \geq 1$ and graph size $n \geq 2k + 1$,
we label the nodes $1, \dots, n$ and connect node $i$ to nodes $i-k, \dots i-1, i+1, \dots, i+k$ modulo $n$ (geometrically, we lay out the $n$ nodes uniformly on a circle and connect each to the $k$ nearest neighbors on either side).
For large graph sizes $n$, the critical quantity for determining global synchronization of $G$ is the density $2k/n$.
Wiley et~al.\ \cite{Wiley2006} show that, if $n,k \to \infty$ with density $2k/n \to \mu \in [0, 1]$,
there is a critical density $\mu_c \approx 0.68$ such that, for large enough $n$, if $\mu < \mu_c$,
there exists an unsynchronized stable equilibrium of \eqref{eq:kuramoto}.
These are regular graphs, so the ordinary and normalized condition numbers are equal: $\frac{\lambda_n(\scrL)}{\lambda_2(\scrL)} = \frac{\lambda_n(L)}{\lambda_2(L)}$.
For this class of graphs, the condition of \Cref{cor:kuramoto_det} is optimal:
\newcommand{\cnfunc}{\kappa}
\begin{theorem}
	\label{thm:circulant}
	Consider the $k$-nearest-neighbor circulant graph $G$ on $n$ edges.
	\begin{enumerate}
		\item \label{it:circ_p1} For all sufficiently large $n$,
		the $k$-nearest-neighbor circulant graph of size $n$ with graph Laplacian $L$ has an unsynchronized stable equilibrium whenever $\frac{\lambda_n(L)}{\lambda_2(L)} > 2$.
		
		\item \label{it:circ_p2} If $n, k \to \infty$ with $\frac{2k}{n} \to \mu \geq 0.6$,
		then $\frac{\lambda_n(L)}{\lambda_2(L)} \to \cnfunc(\mu)$,
		where $\cnfunc$ is an analytic and strictly decreasing function on $[0.6, 1]$ such that $\cnfunc(\mu_c) = 2$.
		
		\item \label{it:circ_p3} Consequently, for every $\epsilon > 0$, there exists a $k$-nearest neighbor circulant graph with an unsynchronized stable equilibrium and $\frac{\lambda_n(L)}{\lambda_2(L)} \leq 2 + \epsilon$.
	\end{enumerate}
\end{theorem}
We provide a proof and further related calculations in \Cref{sec:proofs_circ}.
The key is that, as $A$ and $L$ are circulant matrices, we can calculate their eigenvalues by a discrete Fourier transform.

\subsubsection{\ER{} random graphs and ``random process'' graphs}
\label{sec:kuramoto_er}
One of the main results of Abdalla et~al.\ \cite{Abdalla2026} is that \ER{} random graphs, in the density regime where they are connected with high probability (i.e., the edge probability $p \geq (1 + \epsilon) \frac{\log n}{n}$ for some $\epsilon > 0$), are also (with high probability) globally synchronizing for large size $n$.
\Cref{cor:kuramoto_det}, together with concentration bounds on the normalized graph Laplacian of \ER{} graphs such as that of Hoffman et al.\ \cite{Hoffman2021}, recovers this result (recall that, by a more circuitous route, this is also an implication of our \Cref{thm:z2_bern}).
However, we can say more.
\Cref{cor:kuramoto_det}, combined with a more particular result of \cite{Hoffman2021} (in particular, their Thm.~1.3 and following discussion), recovers the main result of Jain et~al.\ \cite{Jain2026} (which was conjectured by \cite{Abdalla2026}):
that, for large enough $n$, with high probability, the ``random process graph'' (generated by adding one edge at a time, chosen uniformly at random from those missing, starting from the empty graph on $n$ vertices) becomes globally synchronizing as soon as it is connected and remains that way.

\subsubsection{Dense graphs}
Let $G$ be an ordinary (unweighted) graph,
and suppose the minimum vertex degree is $\dmin = \mu(n-1)$ for some density parameter $\mu > 0$.
A long line of work, culminating in \cite{Townsend2020,Kassabov2021,Yoneda2021} (see those papers for further context and references), has sought to understand what is the minimum value of $\mu$ that guarantees global synchrony of $G$.
In particular, Kassabov et~al.\ \cite{Kassabov2021} show that $\mu \geq 0.75$ suffices, and they conjecture that this is optimal.
Evidence for this conjecture is provided by Townsend et~al.\ \cite{Townsend2020},
who construct a sequence of graphs for which $\mu$ approaches $0.75$ from below and the potential function has spurious second-order critical points (although, empirically, these are still \emph{unstable} equilibria of \eqref{eq:kuramoto}).

Our \Cref{cor:kuramoto_det} allows us to recover the result of \cite{Kassabov2021} quite easily.
As the complete graph on $n$ vertices has $\lambda_2(L) = \lambda_n(L) = n$,
monotonicity of $L$ in the graph edges implies that $\lambda_n(L) \leq n$,
and the Gershgorin theorem (noting that each missing edge contributes 2 to the $\ell_1$ norm error of each row of $L$ from that of the complete graph) implies
\begin{align*}
	\lambda_2(L)
	&\geq n - 2 [ (n - 1) - \dmin ] \\
	&= 2 \dmin - n + 2 \\
	&= 2\mu (n-1) - n + 2.
\end{align*}
This lower bound is also given by Fiedler \cite[Property~3.8]{Fiedler1973}.
If
\[
	\mu > \frac{3n - 4}{4 (n-1)} = \frac{3}{4} - \frac{1}{4(n-1)},
\]
we obtain $\lambda_2(L) > \frac{n}{2} \geq \frac{\lambda_n(L)}{2}$.
In particular, this holds if $\mu \geq 0.75$.
Thus we recover the result of \cite{Kassabov2021}.
As our analysis (as well as that in \cite{Kassabov2021}) relies on ruling out spurious second-order critical points of the potential function,
we cannot hope, considering the examples constructed by Townsend et~al.\ \cite{Townsend2020}, that our current methods will give a better density condition than $\mu \geq 0.75$ for large $n$.

\subsubsection{Expander, Ramanujan, and random regular graphs}
\label{sec:expander}
Following Abdalla et~al.\ \cite{Abdalla2026}, we call a graph $G$ on $n$ vertices an $(n, d, \alpha)$-expander for some $d > 0$, $\alpha \geq 0$ if the adjacency matrix $A$ satisfies $\opnorm{A - (d/n) \ones \ones^T} \leq \alpha d$.
One of their results \cite[Thm.~1.4]{Abdalla2026} is that if $G$ is $d$-regular and an $(n,d,\alpha)$-expander for $\alpha \leq 0.0816$, then $G$ is globally synchronizing.
We can improve this result as an easy consequence of \Cref{cor:kuramoto_det}:
\begin{corollary}
	\label{cor:expander_reg}
	If a graph $G$ is $d$-regular and an $(n,d,\alpha)$-expander for $\alpha < 1/3$,
	then $G$ is globally synchronizing.
\end{corollary}
We could also obtain a variant of their result \cite[Thm.~1.10]{Abdalla2026} for non-regular graphs from \Cref{cor:kuramoto_det} and our normalized Laplacian condition number bound \Cref{thm:normalized_cn}, but for brevity we do not pursue this here.

Abdalla et~al.\ \cite{Abdalla2026} furthermore showed, as a consequence of their result for regular expander graphs, that Ramanujan graphs (regular graphs with optimal spectral expander properties) and (with high probability) sufficiently large random regular graphs are globally synchronizing when the vertex degree $d$ is at least $600$.
\Cref{cor:kuramoto_det} allows us to improve this requirement to $d \geq 35$.
Indeed, for any $\epsilon > 0$, $d$-regular Ramanujan graphs and (with high probability) sufficiently large $d$-regular random graphs (see \cite{Friedman2003}) satisfy
\[
	\opnorm*{A - \frac{d}{n} \ones \ones^T} \leq 2 \sqrt{d - 1} + \epsilon.
\]
For small enough $\epsilon$, we have $\frac{2 \sqrt{d - 1} + \epsilon}{d} < 1/3$ whenever $d \geq 35$.
This makes progress toward but is still far from solving another conjecture of Bandeira \cite[Conj.~8]{Bandeira2024oberwolfach}
that random $3$-regular graphs synchronize.

\subsubsection{Random geometric graphs}
Another class of random graph considered in recent papers on oscillator synchronization is the \emph{random geometric graph} (RGG).
Although precise definitions vary, the general model is that each graph vertex has an associated random spatial position in some metric space,
and two vertices are connected by an edge if they are sufficiently close in the distance metric.
See, for example, Abdalla et~al.\ \cite{Abdalla2024} for further background and motivation for oscillator networks on RGGs.
That work \cite{Abdalla2024} has a positive result showing that RGGs on spheres synchronize when they are sufficiently dense (with the average vertex degree growing polylogarithmically in the problem size $n$) and the sphere dimension is large enough (also growing polylogarithmically in $n$).
This analysis is based on a result of Liu et~al.\ \cite{Liu2022} which says that in these parameter regimes the RGG can be well approximated by an \ER{} graph.
The limiting factor appears to be the graph approximation result \cite{Liu2022}, so our \Cref{cor:kuramoto_det} does not likely bring much improvement.
Indeed, in fixed spatial dimension, there are results (e.g., \cite{Silva2011,Adhikari2022}) showing that the (normalized) graph Laplacian does \emph{not} have a small condition number unless the connection range is large (much larger than what is necessary for the graph to be connected).
Thus our techniques are unlikely to yield interesting positive synchronization results for sparse RGGs in a small spatial dimension.
Another negative result is that of De~Vita et~al.\ \cite{DeVita2025},
who prove that the RGG on the circle for connected but sparse parameter regimes is \emph{not} globally synchronizing; qualitatively, the demonstrated stable equilibria are similar to the stable twisted states of nearest-neighbor circulant graphs \cite{Wiley2006}.

\section{A bound on the normalized Laplacian condition number}
\label{sec:LD_condnum}
For many of our applications,
we will apply \Cref{thm:sync_cn} with the ``natural'' diagonal preconditioner
\[
	D = \ddiag(C z z^T) \quad \Longrightarrow \quad L_D = D^{-1/2} L D^{-1/2} = I_n - D^{-1/2} C D^{-1/2}.
\]
In the case that $z = \ones$ and $C$ is the adjacency matrix of a (possibly weighted) graph $G$,
$L_D$ is then the symmetric normalized graph Laplacian.

For the $\Z_2$ synchronization models that we consider in \Cref{sec:intro_appls},
the cost matrix $C$ can be approximated (at least in expectation) by a matrix of the form $\frac{\dbr}{n} z z^T$ for some $\dbr > 0$ (if $z = \ones$, this is, modulo diagonal entries, a rescaled complete-graph adjacency matrix; $\dbr$ can then be interpreted as an average vertex degree).
The normalized Laplacian of this rank-one matrix is simply $I_n - \frac{1}{n} z z^T$, which has a condition number of $1$.

More generally, we consider the case where $C$ is well approximated in a spectral sense by \emph{any} rank-one matrix with the same sign pattern as $z z^T$.\footnote{The fact that rank-one matrices are suitable for synchronization was already noted by Criscitiello et~al.\ \cite{Criscitiello2026} with quite different methods.}
Specifically, for some vector $a \in \R^n$ with strictly positive entries,
let 
\[
	\Cbr = \diag(z) a a^T \diag(z), \quad \text{and} \quad
	\Dbr \coloneqq \ddiag(\Cbr z z^T) = \norm{a}_1 \diag(a).
\]
The corresponding normalized Laplacian matrix is
\[
	\Lbr \coloneqq I_n - \Dbr^{-1/2} \Cbr \Dbr^{-1/2} = I_n - \frac{1}{\norm{a}_1} \diag(z) a^{1/2} (a^{1/2})^T \diag(z),
\]
where $a^{1/2} \in \R^n$ is the elementwise square root of $a$: $a^{1/2}_i = \sqrt{a_i}$.
Crucially, $\lambda_2(\Lbr) = \lambda_n(\Lbr) = 1$.

The following spectral approximation bound for the true normalized Laplacian $L_D$ is the main result of this section.
Qualitatively similar results for random graphs were shown by Chung and Radcliffe \cite{Chung2011} (whose proofs inspired ours), but those results are not quite sufficient for our applications.
\newcommand{\minratio}{\kappa_D}
\newcommand{\cnerrorbd}{\delta_C}
\begin{theorem}
	\label{thm:normalized_cn}
	Let $D = \ddiag(C z z^T)$, and assume $\dmin \coloneqq \min_i~D_{ii} > 0$. Set
	\[
	\minratio \coloneqq \max\{ \opnorm{\Dbr D^{-1}}, 1 \} = \max\braces*{ \parens*{\max_i~\frac{\Dbr_{ii}}{D_{ii}}}, 1}.
	\]
	Then
	\[
		\opnorm{L_D - \Lbr} \leq 2\minratio^2 \opnorm{\Dbr^{-1/2}(C - \Cbr)\Dbr^{-1/2}} \eqqcolon \cnerrorbd.
	\]
	Consequently, if $\cnerrorbd < 1$, $\lambda_2(L_D) > 0$, and
	\[
		\frac{\lambda_n(L_D)}{\lambda_2(L_D)} \leq \frac{1 + \cnerrorbd}{1 - \cnerrorbd}.
	\]
	In the case $\Cbr = \frac{\dbr}{n} z z^T$, we have
	\[
		\cnerrorbd = \frac{2\dbr}{(\dmin \wedge \dbr)^2} \opnorm{C - \Cbr}.
	\]
\end{theorem}
The bound on $\opnorm{L_D - \Lbr}$ is clearly suboptimal by a factor of 2 when $D = \Dbr$ (e.g., $C$ is the adjacency matrix of a $\dbr$-regular graph).
\begin{proof}
	For notational simplicity, we can reduce to the case $z = \ones$ by multiplying $C$ and $\Cbr$ on the left and right by the orthogonal matrix $\diag(z)$, which does not change the spectral properties we are interested in (see also the proofs of \Cref{thm:sync_cn,thm:sync_cn_cplx} in \Cref{sec:proof_condnum}).
	From now on, assume this is the case;
	thus $\Cbr = a a^T$.
	
	We decompose
	\begin{align*}
		L_D - \Lbr &= \Dbr^{-1/2} \Cbr \Dbr^{-1/2} - D^{-1/2} C D^{-1/2} \\
		&= (\Dbr^{-1/2} \Cbr \Dbr^{-1/2} - D^{-1/2} \Cbr D^{-1/2}) + D^{-1/2} (\Cbr - C) D^{-1/2}.
	\end{align*}
	We can easily bound
	\begin{equation}
		\label{eq:opnorm_bd_1}
		\opnorm{ D^{-1/2} (\Cbr - C) D^{-1/2} } \leq \minratio \opnorm{\Dbr^{-1/2} (C - \Cbr) \Dbr^{-1/2}}.
	\end{equation}
	Next, we bound
	\begin{align}
		&\negqquad\opnorm{\Dbr^{-1/2} \Cbr \Dbr^{-1/2} - D^{-1/2} \Cbr D^{-1/2}} \nonumber \\
		&= \opnorm{\Dbr^{-1/2} a a^T \Dbr^{-1/2} - D^{-1/2} a a^T D^{-1/2}} \nonumber \\
		&= \opnorm{ (\Dbr^{-1/2} a - D^{-1/2} a) (\Dbr^{-1/2} a)^T + D^{-1/2} a (\Dbr^{-1/2} a - D^{-1/2} a)^T} \nonumber \\
		&\leq (\norm{\Dbr^{-1/2} a} + \norm{D^{-1/2} a}) \norm{D^{-1/2} a - \Dbr^{-1/2} a} \nonumber \\
		&\leq \parens*{ 1 + \max_i~\sqrt{\frac{\Dbr_{ii}}{D_{ii}} }} \norm{\Dbr^{-1/2} a} \norm{D^{-1/2} a - \Dbr^{-1/2} a } \nonumber \\
		&\leq 2\minratio^{1/2} \norm{D^{-1/2} a - \Dbr^{-1/2} a } \label{eq:opnorm_bd_2},
	\end{align}
	where the last inequality uses the fact that $\norm{\Dbr^{-1/2} a} = 1$.
	
	Finally, we bound
	\begin{align*}
		\norm{D^{-1/2} a - \Dbr^{-1/2} a }^2
		&= \sum_{i=1}^n a_i^2 \parens*{ \frac{1}{\sqrt{D_{ii}}} - \frac{1}{
				\sqrt{\Dbr_{ii}} } }^2 \\
		&\overset{\mathclap{\text{(a)}}}{\leq} \frac{1}{4} \sum_{i=1}^n \frac{a_i^2}{(D_{ii} \wedge \Dbr_{ii})^3} \parens{ D_{ii} - \Dbr_{ii} }^2 \\
		&\leq \frac{\minratio^3}{4 \norm{a}_1^2} \sum_{i=1}^n \frac{(D_{ii} - \Dbr_{ii} )^2}{\Dbr_{ii}} \\
		&= \frac{\minratio^3}{4 \norm{a}_1^2} \norm{\Dbr^{-1/2}(C - \Cbr) \ones}^2 \\
		&= \frac{\minratio^3}{4} \norm{\Dbr^{-1/2}(C - \Cbr) \Dbr^{-1/2} (\Dbr^{-1/2} a)}^2 \\
		&\leq \frac{\minratio^3}{4} \opnorm{\Dbr^{-1/2}(C - \Cbr) \Dbr^{-1/2}}^2.
	\end{align*}
	Inequality (a) comes from calculating the maximum absolute value of the derivative of $x \mapsto \frac{1}{\sqrt{x}}$ between $\Dbr_{ii}$ and $D_{ii}$.
	We used several times the fact that $\Dbr_{ii} = \norm{a}_1 a_i$.
	The last inequality uses $\norm{\Dbr^{-1/2} a} = 1$.
	
	Combined with \eqref{eq:opnorm_bd_2}, we obtain
	\[
	\opnorm{\Dbr^{-1/2} \Cbr \Dbr^{-1/2} - D^{-1/2} \Cbr D^{-1/2}}
	\leq \minratio^2 \opnorm{\Dbr^{-1/2}(C - \Cbr) \Dbr^{-1/2}}.
	\]
	Combining this with \eqref{eq:opnorm_bd_1} and noting that $\minratio \geq 1$ by definition gives the final bound on $\opnorm{L_D - \Lbr}$.
	The condition number bound comes from Weyl's inequality and the fact that $\lambda_2(\Lbr) = \lambda_n(\Lbr) = 1$.
\end{proof}

\section{Proofs for specific applications}
\label{sec:proofs_appls}
\subsection{Random $\Z_2$ synchronization models}
\label{sec:proofs_z2}
In this section, we provide proofs for the statistical $\Z_2$ synchronization results \Cref{thm:z2_gauss,thm:z2_bern,thm:sbm}.
We will use \Cref{thm:sync_cn,thm:normalized_cn};
in particular, it will suffice to show that the quantity $\cnerrorbd$ from \Cref{thm:normalized_cn} is (asymptotically, with high probability) strictly less than $1/3$.
This guarantees that the appropriate Laplacian condition number is strictly less than $2$, with which we can apply \Cref{thm:sync_cn} for all $\sphdim \geq 2$.

To show that, for each problem, $\cnerrorbd < 1/3$,
we use probabilistic concentration inequalities already present in \cite{McRae2025}.
Thus we will simply refer to the appropriate arguments in that paper, discussing any (minor) modifications needed.

We begin with the Gaussian noise model described in \Cref{sec:intro_gaussian}.
\begin{proof}[Proof of \Cref{thm:z2_gauss}]
	Although our proof resembles that of \cite[Cor.~1]{McRae2025}, the concentration inequalities we need are already present in \cite{Bandeira2018}.
	Without loss of generality, assume, for simplicty, that $z = \ones$ (we can reduce to this case by a change of variable: see, e.g., \cite{Bandeira2016a,McRae2024}).
	In the language of \Cref{sec:LD_condnum}, we take $\Cbr = \ones \ones^T$, $\dbr = n$, $\Dbr = n I_n$, and $\Lbr = I_n - \frac{1}{n} \ones \ones^T$.
	
	Standard Gaussian matrix concentration inequalities (see, e.g., \cite[Section 3.2]{Bandeira2018}) imply that, with probability $\to 1$ as $n \to \infty$,
	\begin{align*}
		\opnorm{C - \Cbr}
		= \sigma \opnorm{W}
		\lesssim \sigma\sqrt{n}.
	\end{align*}
	Furthermore, standard deviation bounds on Gaussian random variables and a union bound (again, see \cite[Section 3.2]{Bandeira2018}) imply that, for any $\epsilon' > 0$, with probability $\to 1$ as $n \to \infty$,
	\begin{align*}
		\dmin &= n + \sigma \min_i~(W \ones)_i \\
		&\geq n - \sigma \sqrt{(2 + \epsilon') n \log n}.
	\end{align*}
	If we take $\epsilon' < \epsilon$,
	then the assumption $\sigma \leq \sqrt{\frac{n}{(2 + \epsilon) \log n}}$ implies
	\[
		\dmin \geq \epsilon'' n
	\]
	for some constant $\epsilon'' > 0$ depending on $\epsilon, \epsilon'$.
	We then obtain
	\[
		\cnerrorbd = \frac{2\dbr}{(\dmin \wedge \dbr)^2} \opnorm{C - \Cbr}
		\lesssim \frac{1}{(\epsilon'')^2 n} \sigma \sqrt{n}
		\lesssim \frac{1}{(\epsilon'')^2 \sqrt{\log n}},
	\]
	which goes to zero as $n \to \infty$.
	Therefore, as $n \to \infty$, with probability $\to 1$, we have $\cnerrorbd < 1/3$ and obtain the result.
\end{proof}

We now turn to the $\Z_2$ synchronization model with an \ER{} graph and Bernoulli noise described in \Cref{sec:intro_bern}.
\begin{proof}[Proof of \Cref{thm:z2_bern}]
	Here, we use elements of the proof of \cite[Thm.~2]{McRae2025}.
	Again, without loss of generality, assume $z = \ones$.
	In the language of \Cref{sec:LD_condnum},
	we take
	\[
		\Cbr \coloneqq \E C = \delta p \ones \ones^T,
	\]
	so $\dbr = n \delta p$.
	
	From the proof of \cite[Thm.~2]{McRae2025} (taking, in the notation of that work, $r \to \infty$)
	we have, with probability $\to 1$ as $n \to \infty$,
	\begin{align*}
		\opnorm{C - \Cbr} &\lesssim \sqrt{n p} = \frac{n \delta p}{\sqrt{n p \delta^2}} \leq \frac{\dbr}{\sqrt{\log n}}, \qquad \text{and} \\
		\dmin &\geq \epsilon' n \delta p = \epsilon' \dbr
	\end{align*}
	for some $\epsilon' > 0$ depending on $\epsilon$.
	We then have
	\[
		\cnerrorbd = \frac{2\dbr}{(\dmin \wedge \dbr)^2} \opnorm{C - \Cbr}
		\lesssim \frac{1}{(\epsilon')^2 \dbr} \cdot \frac{\dbr}{\sqrt{\log n}}
		= \frac{1}{(\epsilon')^2 \sqrt{\log n}}.
	\]
	The rest of the proof is like that of \Cref{thm:z2_gauss} above.
\end{proof}

Finally, we consider the stochastic block model (SBM) described in \Cref{sec:intro_sbm}.
\begin{proof}[Proof of \Cref{thm:sbm}]
	We now draw on the proof of \cite[Thm.~3]{McRae2025}.
	For notational simplicity, set the diagonal elements of the adjacency matrix $A$ to be $p$;
	this only affects the optimization problem \eqref{eq:opt_general} via the calculation of the centering constant $\sbmcenterconst$, and this will have asymptotically negligible effect for large $n$ (as made explicit in the proof of \cite[Thm.~3]{McRae2025}).
	
	With this simplification, one can easily verify that
	\[
		\Cbr \coloneqq \E C = \frac{p - q}{2} z z^T.
	\]
	Thus $\dbr = n \frac{p - q}{2}$.
	
	From the proof of \cite[Thm.~3]{McRae2025} (again taking, in the notation of that work, $r \to \infty$)
	we have, with probability $\to 1$ as $n \to \infty$,
	\begin{align*}
		\opnorm{C - \Cbr} &\lesssim \sqrt{np} \leq \frac{n (p - q)}{\sqrt{\log n}} \lesssim \frac{\dbr}{\sqrt{\log n}}, \qquad \text{and} \\
		\dmin &\geq \epsilon' n \frac{p-q}{2} = \epsilon' \dbr
	\end{align*}
	for some $\epsilon' > 0$ depending on $\epsilon$.
	The rest follows exactly as in the above proof of \Cref{thm:z2_bern}.
\end{proof}

\subsection{$k$-nearest-neighbor circulant graphs}
\label{sec:proofs_circ}
In this section, we calculate the spectral properties of the $k$-nearest-neighbor circulant graphs described in \Cref{sec:kuramoto_circ}, providing along the way a proof of \Cref{thm:circulant}.
It is necessary at several points to find the maxima and minima of functions of an integer variable;
similarly to \cite{Wiley2006},
we leave certain steps to straightforward numerical verification rather than giving analytic proofs.

Let $G$ be a $k$-nearest neighbor circulant graph of size $n$.
The adjacency matrix $A$ is a circulant matrix corresponding to circular convolution with the signal (defined on the integers modulo $n$)
\[
	h_A[i] = \sum_{\substack{j=-k \\ j \neq 0}}^k \delta[i - j],
\]
where $\delta[\cdot]$ is the Kronecker delta function.
The eigenvalues of $A$ are the discrete Fourier transform (DFT) coefficients of $h_A$,
which are, by standard identities,
\[
	H_A[m] = \frac{\sin \parens*{ \frac{(2k+1) \pi}{n} m }}{\sin\parens*{ \frac{\pi}{n} m }} - 1,
\]
with $H_A[0] = 2k$ in the usual limiting sense.
$H_A[0] = \sum_i h_A[i]$ is precisely the vertex degree.
The graph Laplacian $L$ then corresponds to circular convolution with
\[
	h_L[i] = H_A[0] \delta[i] - h_A[i].
\]
The eigenvalues of $L$ are then the DFT coefficients
\begin{align*}
	H_L[m]
	= H_A[0] - H_A[m] = 2k + 1 - \frac{\sin \parens*{ \frac{(2k+1) \pi}{n} m }}{\sin\parens*{ \frac{\pi}{n} m }}.
\end{align*}
$H_L[0] = 0$ corresponds to the zero eigenvalue of $L$ with eigenvector $\ones$.
To calculate the condition number of $L$, we must consider the values of $H_L[m]$ for $m \neq 0 \mod n$.
Furthermore, note that $H_L[m] = H_L[n-m]$,
so, for the purpose of calculating the condition number of $L$, we can restrict ourselves to $1 \leq m \leq n/2$.

One can verify that, for large enough $n$ and any valid value of $k$,
the minimum value of $H_L[m]$ for $1 \leq m \leq n/2$ is $H_L[1]$.
Thus
\[
	\lambda_2(L) = H_L[1] = 2k + 1 - \frac{\sin \parens*{ \frac{(2k+1) \pi}{n} }}{\sin\parens*{ \frac{\pi}{n} }}.
\]
To calculate $\lambda_n(L) = \max_m~H_L[m]$,
we need to know what value(s) of $m$ to consider.

First,
we restrict ourselves to $k \geq 0.2 n$.
The precise ratio is not important for now; it suffices to assume that $k$ is proportional to $n$ while including the ``interesting'' parameter regimes covered by \Cref{thm:circulant}.
We will handle smaller values of $k$ for part~\ref{it:circ_p1} of \Cref{thm:circulant} separately below.

In this case, for $1 \leq m \leq n/2$, $H_L[m] = 2k + O(\frac{n}{m}) = 2k + O(\frac{k}{m})$.
Therefore, the largest value of $H_L[m]$
either is approximately $2k$ or occurs for small (bounded) $m$.
The subsequent calculations will indeed show that there are small values of $m$ for which $H_L[m]$ is significantly larger than $2k$.

It thus suffices to consider bounded values of $m$.
For all $m$ up to some sufficiently large but fixed bound,
\newcommand{\HAideal}{\Hbr_{A,\mu}}
\newcommand{\HLideal}{\Hbr_{L,\mu}}
\newcommand{\HAidealkn}{\Hbr_{A,2k/n}}
\newcommand{\HLidealkn}{\Hbr_{L,2k/n}}
\[
	\frac{H_A[m]}{n} =  \frac{\sin \parens*{\pi \frac{2k}{n} m}}{\pi m} + O\parens{n^{-1}} = \HAidealkn[m] + O\parens{n^{-1}} ,
\]
where
\[
	\HAideal[m] \coloneqq \frac{\sin \parens*{\pi \mu m}}{\pi m}.
\]
This approximation holds uniformly in $k$, even for small values of $k$, which will be useful later.
Next, this implies, for bounded $m$, uniformly in $k$,
\[
	\frac{H_L[m]}{n} = \HLidealkn[m] + O(n^{-1}),
\]
where
\[
	\HLideal[m] \coloneqq \HAideal[0] - \HAideal[m] = \mu - \frac{\sin (\pi \mu m)}{\pi m}.
\]
One can verify numerically that, for $\mu \geq 0.59$,\footnote{In what follows, we consider $\mu \geq 0.6$; the $0.01$ of slack is useful for the limiting/approximation arguments.} the largest value of $\HLideal[m]$ for $m \geq 1$ is $\HLideal[2]$ (and it is strictly larger than any other value).
Thus, for large enough $n$ and $k \geq 0.3n$, the same is true of $H_L$.
Then
\[
	\frac{\lambda_n(L)}{\lambda_2(L)} = \frac{H_L[2]}{H_L[1]}.
\]
Therefore, if $n, k \to \infty$ with $\frac{2k}{n} \to \mu \geq 0.6$, we will have
\[
	\frac{\lambda_n(L)}{\lambda_2(L)} \to \frac{\HLideal[2]}{\HLideal[1]} \eqqcolon \cnfunc(\mu).
\]
$\cnfunc$ is clearly analytic as a ratio of trigonometric polynomials, and one can verify numerically that $\cnfunc$ is strictly decreasing on $[0.6, 1]$.
$\cnfunc(\mu) = 2$ when
\[
	2\parens*{ 1 - \frac{\sin (\pi \mu)}{\pi \mu} } = \frac{2 \HLideal[1]}{\mu} = \frac{\HLideal[2]}{\mu} = 1 - \frac{\sin (2 \pi \mu)}{2 \pi \mu}.
\]
The unique solution for $\mu \in [0.6, 1]$ is precisely the critical value $\mu_c \approx 0.68$ found (by solving the same equation) by Wiley et~al.\ \cite{Wiley2006}.
This proves part~\ref{it:circ_p2} of \Cref{thm:circulant}.

Now, we turn to the proof of part~\ref{it:circ_p1}.
Previously, we considered only $k \geq 0.3n$;
one can check numerically that, for $\mu \leq 0.6$, $\frac{\HLideal[3]}{\HLideal[1]} \geq 2.2$,
so, for sufficiently large $n$ and all $k \leq 0.3n$,
we will have $\frac{\lambda_n(L)}{\lambda_2(L)} \geq \frac{H_L[3]}{H_L[1]} > 2$,
and we know from \cite{Wiley2006} that these graphs have spurious stable equilibria.

Therefore, from now on, we can again restrict ourselves to the case $k \geq 0.3n$.
The argument of \cite{Wiley2006} hinges on the stability of the first ``twisted state'' with angles $\theta_i = \frac{2\pi}{n} i$.
The stability of this state
corresponds to $\lambda_2(\Ltl) > 0$,
where $\Ltl$ is the Laplacian matrix%
\footnote{$\Ltl$ is the Hessian matrix of the potential function of \cite{Wiley2006} of which \eqref{eq:kuramoto} is the minimizing gradient flow. The potential discussed in \Cref{sec:intro_kuramoto} is the same modulo a sign change and rescaling. The zero eigenvalue from $\Ltl \ones = 0$ corresponds to a trivial global shift of the angles.}
corresponding to $\Atl$ given by 
\[
	\Atl_{ij} = A_{ij} \cos\parens*{ 2\pi \frac{i - j}{n} }.
\]
$\Atl$ is a circulant matrix corresponding to circular convolution with the signal
\newcommand{\hAtl}{h_{\Atl}}
\newcommand{\HAtl}{H_{\Atl}}
\newcommand{\HLtl}{H_{\Ltl}}
\[
	\hAtl[i] = h_A[i] \cos\parens*{2\pi \frac{i}{n}}.
\]
We then have, by the DFT modulation identity,
\[
	\HAtl[m] = \frac{1}{2}(H_A[m-1] + H_A[m+1]),
\]
which, in turn, implies
\begin{align*}
	\HLtl[m] &= \HAtl[0] - \HAtl[m] \\
	&= H_A[1] - \frac{1}{2}(H_A[m-1] + H_A[m+1]) \\
	&= -H_L[1] + \frac{1}{2}( H_L[m-1] + H_L[m+1] ).
\end{align*}
One can verify numerically (via approximation of $H_L$ by $\HLidealkn$) that, for large enough $n$ and $k \geq 0.3n$, the minimum value of $\HLtl[m]$ for $1 \leq m \leq n/2$ is the first Fourier coefficient
\begin{align*}
	\HLtl[1] = -H_L[1] + \frac{1}{2} H_L[2].
\end{align*}
This agrees with the observation of Wiley et~al.\ \cite{Wiley2006}
that this is the critical Fourier coefficient determining stability.
Recalling from before that the condition number of $L$ is (again, for large enough $n$ and $k \geq 0.3n$) $\frac{\lambda_n(L)}{\lambda_2(L)} = \frac{H_L[2]}{H_L[1]}$,
we see that $\frac{\lambda_n(L)}{\lambda_2(L)} > 2$ implies $\lambda_2(\Ltl) = \HLtl[1] > 0$,
in which case the twisted state is a stable equilibrium.
Thus we obtain part~\ref{it:circ_p1} of \Cref{thm:circulant}.

To prove part~\ref{it:circ_p3},
choose $\mu \in [0.6, \mu_c)$ such that $2 < \cnfunc(\mu) < 2 + \epsilon$;
then, choosing $n,k \to \infty$ such that $\frac{2k}{n} \to \mu$,
for large enough $n$, part~\ref{it:circ_p2} tells us that $2 < \frac{\lambda_n(L)}{\lambda_2(L)} \leq 2 + \epsilon$,
and, from this, part~\ref{it:circ_p1} implies the graphs have spurious stable equilibria (alternatively, we could have used the result of \cite{Wiley2006} directly for this).
This completes the proof of \Cref{thm:circulant}.

\section{Proof of deterministic landscape results}
\label{sec:proof_condnum}
\subsection{Real case}
\label{sec:proof_condnum_real}
In this section, we give a proof of \Cref{thm:sync_cn},
which is an extension of \cite[Thm.~2.2]{Endor2026} to allow for arbitrary diagonal preconditioners.
Although the proof of \rew{} \cite{Endor2026} could likely be adapted to allow for a preconditioner,
we present a new and shorter proof that may be of independent interest.
We make several simplifications compared to \cite{Endor2026},
but the essential ingredients remain the same: the proof combines the condition-number--based argument of Ling \cite{Ling2025} with a novel refinement based on the zero-gradient condition.
We will comment in several places on the parallels and differences among the three proofs.

By the change of variables
\[
	Y \leftarrow \diag(z) Y, \qquad C \leftarrow \diag(z) C \diag(z)
\]
(see, e.g., \cite{Ling2023a,McRae2024} for more details),
we can assume, without loss of generality, that $z = \ones$.
We then denote, for brevity, $L \coloneqq L(\ones) = L(z)$.

From now on, assume $Y$ is a second-order critical point of \eqref{eq:opt_general}.
We want to show that $Y Y^T = z z^T = \ones \ones^T$.

\paragraph{Error quantification}
First, it is useful to decompose the columns of $Y$ into their orthogonal projections onto $\spn\{\ones \}$ and $\spn\{\ones \}^\perp$.
To accommodate the preconditioner $D$,
an innovation in our argument is to understand orthogonality with respect to the inner product $\ipD{x}{y} \coloneqq \ip{x}{D y}$.
Thus we write
\[
	Y = \ones v^T + W, \qquad \text{where} \qquad v \coloneqq \frac{1}{\tr D} Y^T D \ones \in \R^r.
\]
This ensures that $W = Y - \ones v^T$ satisfies
$W^T D \ones = 0$.
Next, note the identity
\begin{equation}
	\label{eq:UUT_decomp}
	\begin{aligned}
		Y Y^T 
		&= \norm{v}^2 \ones \ones^T + W v \ones^T + \ones (Wv)^T + W W^T \\
		&= Y v \ones^T + \ones (Yv)^T + W W^T - \norm{v}^2 \ones \ones^T.
	\end{aligned}
\end{equation}
Furthermore, by the Pythagorean theorem for $\ipD{\cdot}{\cdot}$,
\begin{equation}
	\label{eq:normU_decomp}
	\tr D = \normF{D^{1/2} Y}^2 = (\tr D) \norm{v}^2 + \normF{ D^{1/2} W}^2.
\end{equation}
In particular, note that $\norm{v} \leq 1$.
To show that $Y Y^T = \ones \ones^T$,
it will suffice to show that $\norm{v} = 1$.

We also write
\[
	d \coloneqq D^{1/2} \ones \qquad \text{and} \qquad \PDonepmat \coloneqq I_n - \frac{1}{\tr D} d d^T.
\]
$\PDonepmat$ is the orthogonal projection matrix onto $\spn\{ d \}^\perp$.
Note that $L_D d = 0$, that is, $d$ is the eigenvector of $L_D$ corresponding to $\lambda_1(L_D) = 0$.
Critically, noting that $D^{1/2} Y = d v^T + D^{1/2} W$,
and $(D^{1/2} W)^T d = W^T D \ones = 0$, we have the bound
\begin{equation}
	\label{eq:LUU_ineq}
	\begin{aligned}
		\ip{L}{Y Y^T}
		&= \ip{L_D}{D^{1/2} Y Y^T D^{1/2}}\\
		&\geq \lambda_2(L_D) \ip{\PDonepmat}{D^{1/2} Y Y^T D^{1/2}} \\
		&= \lambda_2(L_D) \normF{D^{1/2} W}^2 \\
		&= \lambda_2(L_D) (\tr D) (1 - \norm{v}^2),
	\end{aligned}
\end{equation}
where the last equality is \eqref{eq:normU_decomp}.

\paragraph{Criticality conditions}
Second-order criticality conditions for \eqref{eq:opt_general} are (see, e.g., \cite{Boumal2019})
\begin{align}
	S(Y) Y &= 0, \quad \text{and} \label{eq:soc_grad} \\
	\ip{S(Y)}{\Ydt \Ydt^T} &\geq 0\ \quad \text{for all} \quad \Ydt \in \tansp \coloneqq\{ \Ydt \in \R^{n \times \sphdim} : \diag(Y \Ydt^T) = 0 \}, \label{eq:soc_hess}
\end{align}
where
\[
	S(Y) \coloneqq \ddiag(C Y Y^T) - C = L - \ddiag(L YY^T).
\]
$\tansp$ is the tangent space to the product of spheres $(\Sdm)^n$ (which is a Riemannian manifold) at the point $Y$.
The condition \eqref{eq:soc_grad} is equivalent to the Riemannian gradient of the objective function on the constraint manifold being zero;
the condition \eqref{eq:soc_hess} is equivalent to the Riemannian Hessian being negative semidefinite.

Note that if $Y Y^T = z z^T = \ones \ones^T$ (which is what we are trying to prove),
we would have $S(Y) = L$, which was precisely the dual certificate matrix $S$ introduced in \Cref{sec:intro_opt}.
However, we cannot assume this.

\paragraph{Using the Hessian condition}
The randomized argument of \cite{Mei2017,Ling2023a,McRae2025,Ling2025,Endor2026}
is to choose, in the Hessian condition \eqref{eq:soc_hess},
\[
	\Ydt = \Ptansp(\ones \gamma^T) = \ones \gamma^T - \ddiag(\ones \gamma^T Y^T)Y,
\]
where $\gamma \in \R^\sphdim$ is a standard normal random vector.
$\Ptansp$ is the (Euclidean) orthogonal projector onto $\tansp$ in $\R^{n \times \sphdim}$ (again, see, e.g., \cite{Boumal2019}).
One can check that
\[
	\E \Ydt \Ydt^T = (\sphdim - 2) \ones \ones^T + (Y Y^T)^{\circ 2}.
\]
Plugging $\Ydt$ into the Hessian condition \eqref{eq:soc_hess} and taking an expectation gives $\ip{S(Y)}{\E \Ydt \Ydt^T} \geq 0$.
Using the facts that the diagonal elements of $\E \Ydt \Ydt^T$ are all equal to $\sphdim-1$ and that $L \ones = 0$, we obtain
\begin{equation}
	\label{eq:ineq_soc}
	(\sphdim-1) \ip{L}{Y Y^T} = \ip{\ddiag(L Y Y^T)}{\E \Ydt \Ydt^T} \leq \ip{L}{(Y Y^T)^{\circ 2}}.
\end{equation}
We could lower bound the left-hand side by \eqref{eq:LUU_ineq}.
At this point, Ling \cite{Ling2025}, in the case $D = I_n$, upper-bounded the right-hand side of \eqref{eq:ineq_soc} by noting that
\begin{align*}
	\ip{L}{(Y Y^T)^{\circ 2}}
	&\leq \lambda_n(L) \ip{\Ponepmat}{(Y Y^T)^{\circ 2}} \\
	&= \lambda_n(L)\parens*{n - \frac{1}{n} \normF{Y Y^T}^2 } \\
	&\leq 2 \lambda_n(L) \ip{\Ponepmat}{Y Y^T}.
\end{align*}
See \cite[Sec.~3.1]{Ling2025} for more details.
However, this argument is not tight enough to obtain the result of \cite{Endor2026} and \Cref{thm:sync_cn}.

\paragraph{Using the zero-gradient condition}
A key innovation of \rew{} \cite{Endor2026} is to use the zero-gradient condition \eqref{eq:soc_grad} to obtain an improved upper bound on a quantity similar to $\ip{L}{(Y Y^T)^{\circ 2}}$.
The following argument, though different in appearance, retains this essential step.

Considering the decomposition of $Y$ into orthogonal components, $Y = \ones v^T + W$,
and the fact that larger $\norm{v}$ corresponds to less error (see \eqref{eq:normU_decomp}),
it is natural to study a perturbation of $Y$ in the direction $\ones v^T$.
To this end, set
\[
	\Ydt_0 = \Ptansp(\ones v^T) = \ones v^T - \ddiag(\ones v^T Y^T) Y = \ones v^T - \diag(Y v) Y.
\]
Note that a similar quantity appears in \cite[eq.~(21)]{Endor2026} (their $\delta e_1$ corresponds to our $v$).
Then, using $L \ones = 0$, $\ddiag(\Ydt_0 Y^T) = 0$, and $S(Y) Y = 0$ (from \eqref{eq:soc_grad}), we have
\begin{equation}
	\label{eq:eq_foc}
	0 = \ip{S(Y)}{\Ydt_0 Y^T} = - \ip{L}{\diag(Yv) Y Y^T}.
\end{equation}
Adding \eqref{eq:eq_foc} to the right-hand side of \eqref{eq:ineq_soc} (with scaling factor $2$, the reason for which will quickly become clear) and recalling \eqref{eq:UUT_decomp},
we obtain
\begin{align*}
	(\sphdim-1) \ip{L}{YY^T}
	&\leq \ip{L}{(Y Y^T)^{\circ 2}} - 2 \ip{L}{\diag(Yv) Y Y^T} \\
	&= \ip{L}{ (Y Y^T - [ (Yv) \ones^T + \ones (Yv)^T ]) \circ (Y Y^T) } \\
	&= \ip{L}{ [ W W^T - \norm{v}^2 \ones \ones^T  ] \circ Y Y^T  }.
\end{align*}
This implies
\begin{equation}
	\label{eq:ineq_basic}
	(\sphdim - 1 + \norm{v}^2) \ip{L}{Y Y^T} \leq \ip{L}{ \underbrace{(W W^T) \circ (Y Y^T)}_{\eqqcolon M}  }.
\end{equation}
The matrix $M$ corresponds (up to an orthogonal projection) to the matrix called ``$M_*$'' in \cite{Endor2026},
although our derivation is quite different.
They write $M_*$ as an arbitrary linear combination of several of the above matrices and then search carefully for the best coefficients;
in our argument, the coefficients follow naturally from the decomposition \eqref{eq:UUT_decomp}.

\paragraph{Final calculations}
The left-hand side of \eqref{eq:ineq_basic} can be lower bounded by \eqref{eq:LUU_ineq}.
The right-hand side can be upper bounded (noting that $M \succeq 0$ by the Schur product theorem) as
\begin{equation}
	\label{eq:psd_ub}
	\begin{aligned}
		\ip{L}{ M }
		&= \ip{L_D}{D^{1/2} M D^{1/2}} \\
		&\leq \lambda_n(L_D) \ip{ \PDonepmat }{ D^{1/2}M D^{1/2} } \\
		&= \lambda_n(L_D) \parens*{ \tr (D^{1/2} M D^{1/2}) - \frac{1}{\tr D} \ip{d d^T}{D^{1/2} M D^{1/2}} },
	\end{aligned}
\end{equation}
We can easily calculate, by \eqref{eq:normU_decomp},
\begin{equation}
	\label{eq:trM}
	\tr (D^{1/2} M D^{1/2}) = \normF{D^{1/2} W}^2 = (\tr D) (1 - \norm{v}^2).
\end{equation}
Finally, noting that $W^T D Y = W^T D W$, we can bound
\begin{align*}
	\ip{d d^T}{D^{1/2} M D^{1/2}}
	&= \ip{D^{1/2} W W^T D^{1/2}}{ D^{1/2} Y Y^T D^{1/2}} \\
	&= \normF{D^{1/2} W W^T D^{1/2}}^2 \\
	&\geq \frac{1}{\sphdim} \tr^2( D^{1/2} W W^T D^{1/2} )\\
	&= \frac{1}{\sphdim} \normF{D^{1/2} W}^4 \\
	&= \frac{(\tr D)^2 }{\sphdim} (1 - \norm{v}^2)^2.
\end{align*}
Our use of a simple trace/Frobenius norm inequality here is quite different from the calculations in \cite{Endor2026}.
Combining this with \eqref{eq:psd_ub} and \eqref{eq:trM} gives
\begin{align*}
	\ip{L}{M} &\leq \lambda_n(L_D) (\tr D) \parens*{ 1 - \norm{v}^2 - \frac{1}{\sphdim} (1 - \norm{v}^2)^2 } \\
	&= \frac{\lambda_n(L_D)}{\sphdim} (\tr D) (\sphdim - 1 + \norm{v}^2)(1 - \norm{v}^2).
\end{align*}
Combining this with \eqref{eq:ineq_basic} and \eqref{eq:LUU_ineq} and dividing by $\tr D$, we obtain
\[
(\sphdim - 1 + \norm{v}^2) \lambda_2(L_D) (1 - \norm{v}^2) \leq \frac{\lambda_n(L_D)}{\sphdim} (\sphdim - 1 + \norm{v}^2)(1 - \norm{v}^2).
\]
Clearly, when $\lambda_2(L_D) > \frac{\lambda_n(L_D)}{\sphdim}$, the only way this is satisfied is if $\norm{v} = 1$.
This implies $Y Y^T = \ones \ones^T$,
which completes the proof.

\subsection{Complex case}
\label{sec:proof_condnum_cplx}
In this section, we provide a proof of \Cref{thm:sync_cn_cplx}.
The proof is similar to that of \Cref{thm:sync_cn} (the real case) in the previous section with a few additional subtleties.
See \cite{Pumir2018} for details on how to extend the Riemannian geometry calculations to the complex case.
We use the convention that the (matrix) inner product is always real,
that is, $\ip{A}{B} = \real(\tr(A^* B))$.

Similarly to the real case, the change of variable
\[
	Y \leftarrow \diag(z)^* Y, \qquad C \leftarrow \diag(z)^* C \diag(z)
\]
allows us to assume, without loss of generality, that $z = \ones$,
and we again denote $L \coloneqq L(\ones) = L(z)$.
Second-order optimality conditions for \eqref{eq:opt_cplx} have a nearly identical form as before:
\begin{align}
	S(Y) Y &= 0, \quad \text{and} \label{eq:grad_cond_cplx} \\
	\ip{S(Y)}{\Ydt \Ydt^*} &\geq 0 \quad\ \text{for all} \quad  \Ydt \in \tansp \coloneqq\{ \Ydt \in \C^{n \times \sphdim} : \real(\diag(Y \Ydt^*)) = 0 \}, \label{eq:hess_cond_cplx}
\end{align}
where
\[
	S(Y) \coloneqq \real(\ddiag(C Y Y^*)) - C = L - \real(\ddiag(L Y Y^*)).
\]

To apply the Hessian condition \eqref{eq:hess_cond_cplx}, we plug in, similarly to the real case, the random tangent vector
\[
	\Ydt = \Ptansp(\ones \gamma^*) = \ones \gamma^* - \real(\ddiag(\ones \gamma^* Y^*)) Y
	= \ones \gamma^* - \diag(\real(Y \gamma)) Y,
\]
where $\gamma$ is now a \emph{complex} standard normal random vector in $\C^\sphdim$.
Writing $Y$ rowwise as
\[
	Y = \begin{bmatrix*} y_1^* \\ \vdots \\ y_n^* \end{bmatrix*},
\]
where $y_1, \dots, y_n \in \C^\sphdim$,
we can calculate
\begin{align*}
	\E_\gamma (\Ydt \Ydt^*)_{ij}
	&= \E_\gamma [ (\gamma^* - \real(y_i^* \gamma) y_i^*)(\gamma^* - \real(y_j^* \gamma)y_j^*)^*  ] \\
	&= \E_\gamma\brackets*{ \parens*{ \gamma^* - \frac{1}{2}(y_i^* \gamma + \gamma^* y_i)y_i^* } \parens*{ \gamma^* - \frac{1}{2}(y_j^* \gamma + \gamma^* y_j)y_j^* }^*  } \\
	&= \E_\gamma \Bigg[ \norm{\gamma}^2 - \frac{1}{2} \parens*{ \gamma^* y_j \gamma^* y_j + \gamma^* y_j y_j^* \gamma + y_i^* \gamma y_i^* \gamma + \gamma^* y_i y_i^* \gamma } \\
	&\qquad \qquad + \frac{1}{4} \parens*{ y_i^* \gamma y_i^* y_j \gamma^* y_j + \gamma^* y_i y_i^* y_j y_j^* \gamma + y_i^* \gamma y_i^* y_j y_j^* \gamma +  \gamma^* y_i y_i^* y_j \gamma^* y_j } \Bigg] \\
	&= \sphdim - \frac{1}{2} \parens*{ \norm{y_i}^2 + \norm{y_i}^2 } + \frac{1}{4} \parens*{ (y_i^* y_j)^2 + \abs{y_i^* y_j}^2 } \\
	&= \sphdim - 1 + \frac{1}{2} \real(y_i^* y_j) y_i^* y_j.
\end{align*}
Many terms in the long middle expression disappear due to the circular symmetry of entries of $\gamma$; any term involving $\gamma$ twice or $\gamma^*$ twice has zero expectation.
We could also obtain the same result as a special case of a calculation in \cite[Sec.~5]{McRae2024}.
Noting that $(Y Y^*)_{ij} = y_i^* y_j$, we have
\[
\E_\gamma \Ydt \Ydt^*
= (\sphdim-1) \ones \ones^* + \frac{1}{2} \real(Y Y^*) \circ (Y Y^*).
\]
Plugging this into \eqref{eq:hess_cond_cplx} with an expectation,
we obtain (rescaling by $2$)
\begin{equation}
	\label{eq:ineq_soc_cplx}
	(2r-1) \ip{L}{Y Y^*} \leq \ip{L}{\real(Y Y^*) \circ (Y Y^*)}.
\end{equation}
To use the zero-gradient equality \eqref{eq:grad_cond_cplx},
we decompose, similarly to the real case,
\[
	Y = \ones v^* + W, \qquad \text{where} \qquad v \coloneqq \frac{1}{\tr D} Y^* D \ones \in \C^r,
\]
which ensures that $W^* D \ones = 0$.
We use the tangent vector
\[
\Ydt_0 \coloneqq \Ptansp(\ones v^*) = \ones v^* - \real(\ddiag(\ones v^* Y^*)) Y
= \ones v^* - \diag(\real(Y v)) Y.
\]
By \eqref{eq:grad_cond_cplx} and the facts that $L \ones = 0$ and $\real(\ddiag(\Ydt_0 Y^*)) = 0$,
\begin{equation}
	\label{eq:eq_foc_cplx}
	0 = \ip{S(Y)}{\Ydt_0 Y^*} = - \ip{L}{\diag(\real(Yv)) Y Y^*}.
\end{equation}
Noting (similarly to the real case) that
\[
Y Y^* = Y v \ones^* + \ones (Yv)^* + W W^* - \norm{v}^2 \ones \ones^*,
\]
combining \eqref{eq:ineq_soc_cplx} and \eqref{eq:eq_foc_cplx} (scaled by $2$) gives
\begin{align*}
	(2r-1)\ip{L}{Y Y^*}
	&\leq \ip{L}{\real(Y Y^*) \circ (Y Y^*)} - 2  \ip{L}{\diag(\real(Yv)) Y Y^*} \\
	&= \ip{L}{ \real( Y Y^* - (Yv \ones^* + \ones (Yv)^*) ) \circ (YY^*) } \\
	&= \ip{L}{ \real( W W^* - \norm{v}^2 \ones \ones^* ) \circ (YY^*) },
\end{align*}
from which we obtain the key inequality
\begin{equation}
	\label{eq:ineq_basic_cplx}
	(2 \sphdim - 1 + \norm{v}^2) \ip{L}{Y Y^*} \leq \ip{L}{ \underbrace{\real(W W^*) \circ (Y Y^*)}_{\eqqcolon M}  }.
\end{equation}
The remainder is almost the same as in the real case, so we omit some details and definitions that are the same as before.
We have, again,
\[
	\tr (D^{1/2} M D^{1/2}) = \normF{D^{1/2} W}^2 = (\tr D) (1 - \norm{v}^2).
\]
Next, noting that $\real(W W^*) D \ones = \real(W W^* D \ones) = 0$,
\begin{align*}
	\ip{D^{1/2} M D^{1/2}}{d d^*}
	&= \ip{D^{1/2} \real(W W^*) D^{1/2}}{D^{1/2} Y Y^* D^{1/2}} \\
	&= \ip{D^{1/2} \real(W W^*) D^{1/2}}{D^{1/2} W W^* D^{1/2}} \\
	&= \normF{D^{1/2} \real(W W^*) D^{1/2}}^2 \\
	&\geq \frac{1}{2r} \tr^2( D^{1/2} \real(W W^*) D^{1/2} )\\
	&=\frac{1}{2r} \normF{D^{1/2} W}^4 \\
	&= \frac{(\tr D)^2}{2r} (1 - \norm{v}^2)^2.
\end{align*}
A subtlety in the trace/Frobenius norm inequality is that $\real(W W^*)$ can have rank up to $2r$.
Finally, similarly to the real case, we obtain the inequality
\begin{align*}
	\ip{L}{M} &\leq \lambda_n(L_D) (\tr D)\parens*{ 1 - \norm{v}^2 - \frac{(1 - \norm{v}^2)^2}{2r} } \\
	&= \frac{\lambda_n(L_D)}{2r} (\tr D)(2r - 1 + \norm{v}^2)(1 - \norm{v}^2).
\end{align*}
Via \eqref{eq:ineq_basic_cplx}, the rest of the proof follows exactly as in the real case.

\section*{Acknowledgments}
The author thanks Pedro Abdalla, Afonso Bandeira, Nicolas Boumal, Faniriana Rakoto Endor, and Irène Waldspurger 
for helpful discussions and inspiration.
The author also thanks the two anonymous referees whose comments led to many improvements in the paper.

\bibliographystyle{siamplain}
\bibliography{./refs}

\end{document}